\documentclass[11pt,a4paper, twoside, openright]{book}
\pdfoutput=1
 \usepackage[T1]{fontenc}
\usepackage{bm}
\usepackage{amsmath,amssymb}
\usepackage{cite}
\usepackage{graphicx}
\usepackage[francais]{babel}
\usepackage[latin1]{inputenc}
\numberwithin{figure}{section}
\numberwithin{equation}{section}
\numberwithin{table}{section}

\usepackage[Lenny]{fncychap}
\usepackage{courier}
\usepackage{type1cm}         

\usepackage{makeidx}         
 \usepackage{graphicx, psfrag}
\usepackage{multicol}        
\usepackage[bottom]{footmisc}

  \newcommand{\ds}{\displaystyle}

\makeindex             
  \newcommand{\be}{\begin{equation}}
\newcommand{\ee}{\end{equation}}

\newcommand{\ba}{\begin{eqnarray}}
\newcommand{\ea}{\end{eqnarray}}
\newcommand{\bal}{\begin{align}}
\newcommand{\eal}{\end{align}}
 
  \newcommand{\un}{\underline}
	\newcommand{\ov}{\overline} 
	\def\BbR{{\mathbb R}}

	\usepackage{amsfonts}				
\newcommand{\lb}{\label}
\def\BbC{\mathbb C}
\def\BbS{\mathbb S}
\def\BbT{\mathbb T}

\def\m{\mathcal}
\newcommand{\fii}{\varphi}
\newcommand{\lm}{\lambda}
\newcommand{\non}{\nonumber}

\begin{document}


\begin{titlepage}
\begin{center}
\Large
\textsc{Ecole Sup\'erieure Priv\'ee d'A\'eronautique et des Technologies de Tunis }\\

\vspace{5cm}


\LARGE
\textsc{\textbf{TRAVAUX PRATIQUES - EXERCICES ET PROBLEMES DE GEODESIE}}\\[0.5\baselineskip]
Collect\'es Par\\[0.5\baselineskip]
\textbf{Abdelmajid BEN HADJ SALEM}\\
\vspace{0.5cm}
\author{\Large \textsc{Ing\'enieur G\'eographe G\'en\'eral }\\
\normalsize
\vspace{1cm}
\textsc{D\'ecembre 2016}
\\ 
\vspace{1cm}
\textsc{Version 3.a}\\
\vspace{4cm}
\textit{abenhadjsalem@gmail.com}}

\end{center}
\end{titlepage}
\tableofcontents 
\chapter*{Pr\'eface}

\addcontentsline{toc}{section}{Pr\'eface}
              Le pr\'esent document est une collection de travaux pratiques, d'exercices et de probl\`emes de g\'eod\'esie pour les \'etudiants en g\'eod\'esie cycle des ing\'enieurs. 
												\newpage
\part{\textsc{Travaux Pratiques de G\'eod\'esie}}
\chapter{Courbes et Surfaces - Ellipse et Ellipsoide de r\'evolution}
\section{\textsc{Courbes et Surfaces}}
\subsection{Courbes Gauches}
* Le tri\`edre de Fren\^et.

* D\'efinition de la courbure.

* Calcul dans le cas d'une courbe plane.
\subsection{Surfaces param\'etr\'ees}
* Plan tangent, vecteur unitaire normal.

* 1\`ere forme quadratique fondamentale, courbes orthogonales sym\'etriques.

* Courbe trac\'ee sur une surface, tri\`edre de Darboux.

* Th\'eor\`eme de Meusnier.

* Directions et courbures principales, courbure d'une section normale.

* Courbures principales d'une surface de r\'evolution.
\section{\textsc{L'Ellipse et l'ellipsoide de r\'evolution}}
* D\'efinition.

* D\'efinitions concernant l'aplatissement, les excentricit\'es.

* D\'efinitions des latitudes:

- param\'etrique : $\psi$,

- g\'eographique: $\varphi$,

- g\'eocentrique: $\omega$,

\subsection{Exercices}
1.  Calculer les composantes du vecteur normal ext\'erieur \`a l'ellipsoide, en d\'eduire les relations (dans les deux sens) entre les lignes trigonom\'etriques de $\varphi$ et celles de $\psi$.

2. Donner les \'equations param\'etriques de l'ellipse et de l'ellipsoide en fonction, repectivement, de $\varphi$ et de $\lambda$ et $\varphi$.

3. Etablir une relation diff\'erentielle entre $\psi$ et $\varphi$.

4. Calculer la diff\'erentielle $d\beta$ de l'arc d'ellipse en fonction de $\varphi$, puis la premi\`ere forme quadratique de l'ellipsoide.

5. Calculer les courbures principales de l'ellipsoide de r\'evolution.

6. Trouver la coordonn\'ee curviligne de l'ellipsoide de r\'evolution qui forme avec la longitude un couple de coordonn\'ees sym\'etriques et qui s'annulle  le long de l'\'equateur.

\section{\textsc{Calcul d'un arc d'ellipse}}
\subsection{Int\'egrales de Wallis}
Soit \`a calculer:
$$ W_{2p}=\int_0^{\Omega}sin^{2p}\omega d \omega $$
On pose:
$$ I_{p-2  }(\Omega)=\int_0^{\Omega}sin^{p-2}\omega cos^2\omega d \omega $$
1. Etablir les formules suivantes:
\ba
W_{p}= W_{p-2}-I_{p-2} \\
(p-1)I_{p-2}=sin^{p-1}\Omega cos \Omega +W_p \\
W_{p}= \ds \frac{p-1}{p}W_{p-2}-\frac{1}{p}sin^{p-1}\Omega cos \Omega
\ea
2. Pr\'eciser la valeur de $W_0$, et proposer un programme (en Matlab) de calcul de $W_{2p}$.
\subsection{Calcul d'un arc d'ellipse m\'eridienne}
On a:
$$ \beta(\varphi)=\int_0^{\varphi}\rho d\varphi $$
avec :
$$ \rho=\frac{a(1-e^2)}{w^3},\quad w^2=1-e^2sin^2\varphi $$
1. D\'evelopper $w^{-3}$ suivant les puissances croissantes de $esin\varphi$.

2. Calculer $\beta(\varphi)$ en fonction des $W_{2p}(\varphi)$.

3. Majorer l'erreur de calcul, lorsqu'on arr\^ete le d\'eveloppement au terme $e^{2n}$. Calculer $n$ si l'on recherche la pr\'ecision du millim\`etre sur $\beta$, quelle que soit la latitude $\varphi$ entre $-\frac{\pi}{2}$ et $\frac{\pi}{2}$.

4. Proposer un organigramme de calcul.

5. Envisager la solution du probl\`eme inverse: calcul de $\varphi$ connaissant $\beta$.

\chapter{Passage des coordonn\'ees cart\'esiennes en ellipsoidiques}

	\section{\textsc{Introduction}}
On consid\`ere un r\'ef\'erentiel g\'eod\'esique $(O,OX,OY,OZ)$ avec un ellipsoide de r\'ef\'erence $\m E $ de param\`etres $a,e^2$ o\`u respectivement $a$ le demi-grand axe et le carr\'e de la premi\`ere excentricit\'e. Les coordonn\'ees cart\'esiennes $(X,Y,Z)$ d'un point $M(\fii,\lm,h)$ sont donn\'ees par:
\be
M\left\{\begin{array}{l}
X=(N+h)cos\fii cos\lm \\
Y=(N+h)cos\fii sin\lm \\
Z=(N(1-e^2)+h)sin\fii
\end{array}\right.
\ee
On veut \'etudier le passage de $(X,Y,Z)$ \`a $(\fii,\lm,h)$.

1. Montrer en consid\'erant que $Z\geq 0$, $h$ v\'erifie:
\be
h\geq -(N(1-e^2)
\ee
Le calcul de $\lm$ est facile, et on a:
\be
tg(\lm)=\ds \frac{Y}{X}
\ee
Le calcul de $\fii,h$ est plus complexe. Sa r\'esolution peut se faire par les trois m\'ethodes suivantes:
 
	 a - les algorithmes it\'eratifs,
	
	 b - les algorithmes finis,
	
	 c - les d\'eveloppements limit\'es.
\subsection{Les Algorithmes It\'eratifs}
  Les algorithmes it\'eratifs conduisent \`a r\'esoudre une \'equation de la forme:
	\be
	x=f(x) \lb{a1}
	\ee
	o\`u $x$ est une inconnue auxiliaire. Soit $\overline{x}$ une solution de (\ref{a1}). Partant d'une solution approch\'ee $x_0$ de $\overline{x}$, on calcule successivement:
	\ba
	x_1=f(x_0) \non \\
	x_2=f(x_1) \non \\
		\vdots  \non \\
	x_i=f(x_{i-1}) \non
	\ea
La m\'ethode converge si, pour un ensemble de voisinages $V_i$ de la solution $\ov{ x}$:
\be
V_i=f(V_{i-1})\subset V_{i-1}
\ee
On dit que $f$ est une fonction contractante au voisinage de $\ov{x}$.

Si $f$ est continue et d\'erivable au voisinage de $\ov x$, la condition de convergence est:
\be
\left. \begin{array} {l} 
\exists \,\,\mbox{un voisinage}\,\,V\,\,\mbox{de}\,\,\ov x \\
\\
\exists k\in \left[0,1\right]
 \end{array}\right\}  \mbox{tels que}\,\,x\in V \Longrightarrow \left|f'(x)\right| \leq k 
\ee
En d'autres termes, la m\'ethode it\'erative converge si $|f'(x)|$ est major\'ee par un coefficient $k$ inf\'erieur \`a 1 au voisinage de la solution. La divergence de la m\'ethode n'implique pas l'inexistance de la solution. 
\subsubsection{Etude de la pr\'ecision de la M\'ethode}
de:
\ba
x_1=f(x_0) \non \\
	x_2=f(x_1) \non \\
		\vdots  \non \\
	x_i=f(x_{i-1}) \non
	\ea
	on tire:
	\ba
\ov x - x_1=f(\ov x)- f(x_0)=(\ov x -x_0)f'(\theta_1) \non \\
\ov x - x_2=f(\ov x)- f(x_1)=(\ov x -x_1)f'(\theta_2) \non \\
		\vdots    \\
\ov x - x_{i-1}=f(\ov x)- f(x_{i-2})=(\ov x -x_{i-2})f'(\theta_{i-1}) \non \\
\ov x - x_{i}=f(\ov x)- f(x_{i-1})=(\ov x -x_{i-1})f'(\theta_{i}) \non  
	\ea
avec, \`a chaque fois $\theta_i\in \left[x_i,\ov x\right]$. On en d\'eduit les $i+1$ in\'egalit\'es:
	\ba
|\ov x - x_0|\leq|\ov x -x_0| \non \\
|\ov x - x_1|<k|\ov x -x_0| \non \\
		\vdots    \non \\
|\ov x - x_{i-1}|<k|\ov x -x_{i-2}| \non \\
|\ov x - x_{i}|<k|\ov x -x_{i-1}| \non  
	\ea
Soit par multiplication membre \`a membre:
\be
|\ov x - x_{i}|<k^i|\ov x -x_0|
\ee
Comme la fonction est contractante, on peut, dans cette \'equation, remplacer $|\ov x-x_0|$ par $|x'_0-x_0|$, $x'_0$ \'etant une valeur approch\'ee de $\ov x$ formant, avec $x_0$ un encadrement de $\ov x$: 
$$\ov x \in \left[ x_0,x'_0\right]$$
L'erreur de la m\'ethode est alors major\'ee par:
\be
|\ov x - x_{i}|<k^i|x'_0 -x_0|
\ee
Le nombre $i$ d'it\'erations n\'ecessaires pour calculer $\ov x$ avec une pr\'ecision $\epsilon_x$ donn\'ee \`a l'avance est obtenu par:
$k^i|x'_0 -x_0|< \epsilon_x $ soit, en se souvenant que $k<1$
\be
\fbox{ $i>\ds \frac{Log \ds \frac{\epsilon_x}{|x'_0 -x_0|}}{Logk} $}
\ee
\subsubsection{Etapes de l'Analyse et la programmation d'une M\'ethode it\'erative:}
L'Analyse et la programmation d'une M\'ethode it\'erative peuvent comporter les \'el\'ements suivants:
\begin{itemize}
	\item D\'emonstration de l'\'equation employ\'ee: $x=f(x)$.
	\item Calcul et majoration de $f'$, eventuellement \'etude des cas o\`u $|f'(x)|>1$.
	\item Calcul de $i$, nombre d'it\'erations n\'ecessaires. Ce nombre peut \^etre fix\'e une fois pour toutes, ou calcul\'e par le programme.
	\item Calcul de $\fii$ et de $h$ par une formule peu sensible aux erreurs sur $x$.
	\item R\'edaction d'un organigramme et d'un programme en Matlab, accompagn\'es d'une notice d'emploi sur les limites d'emploi, le temps d'\'execution et l'encombrement en machines. 
	\item Un jeu d'essai au moins.
\end{itemize}
 
\subsubsection{Les M\'ethodes it\'eratives propos\'ees:}
\begin{enumerate}
	\item $\fii=Arctg\left(\ds \frac{Z}{p}+\frac{Ne^2sin\fii}{p}\right)$ avec $p=\sqrt{X^2+Y^2}$.
	\item $\fii=Arctg\left[\ds \frac{Z}{p}\left(1-\frac{Ne^2cos\fii}{p}\right)^{-1}\right]$ avec $p=\sqrt{X^2+Y^2}$.
	\item $\fii=\psi+Arcsin\ds \frac{Ne^2sin2\fii }{2r}$, avec $\psi=Arctg\ds \frac{Z}{\sqrt{X^2+Y^2}},\,r=\sqrt{p^2+Z^2}$.
\end{enumerate}
\subsection{Les Algorithmes finis}
Les algorithmes finis conduisent \`a la r\'esolution d'une \'equation du 4\`eme degr\'e:
\be
x^4+a_1x^3+a_2x^2+a_3x+a_4=0 \lb{b11}
\ee
o\`u $x$ est une variable auxiliaire. On donne ci-dessous quelques indications pour la r\'esolution de l'\'equation (\ref{b11}).

On \'elimine le terme du troisi\`eme degr\'e par la transformation lin\'eaire:
\be
y=x+\ds \frac{a_1}{4}
\ee
L'\'equation (\ref{b11}) devient:
\be
y^4+a'_2y^2+a'_3y+a'_4=0 \lb{b12}
\ee
On abaisse ensuite le degr\'e de cette \'equation en posant:
\be
2y=u+v+w \lb{b3}
\ee
Entre les variables ind\'ependantes $u,v,w$, on peut encore imposer deux relations, par exemple:
\ba
u^2+v^2+w^2=-2a'_2\lb{b5}\\
u.v.w=-a'_3\lb{b6}
\ea
En utilisant les fonctions sym\'etriques des racines, on montre que $u^2,v^2,w^2$ sont solutions de l'\'equation:
\be
z^3+2a'_2z^2+(a'^2_2-4a'_4)z-a'_3=0 \lb{b4}
\ee
On fait dispara\^itre le terme du second degr\'e par le changement:
\be
\xi=z+\ds \frac{2a'_2}{3}\lb{b7}
\ee
et on r\'esout l'\'equation par la m\'ethode de Cardan (Chercher dans le Web). 

L'\'equation (\ref{b4}) admet trois solutions dans le corps des nombres complexes $\BbC$: $u^2,v^2,w^2$. Ce-ci fournit \textit{\`a priori} 8 possibilit\'es de calculer $y$ par (\ref{b3}). Cependant, le signe du produit $uvw$ est impos\'e par (\ref{b6}), et l'\'equation (\ref{b11}) n'admet bien que quatre solutions parmis lesquelles il faudra choisir.
\subsubsection{L'Analyse et la Programmation d'une M\'ethode Finie:}
L'Analyse et la programmation d'une m\'ethode finie comportera les \'el\'ements suivants:
\begin{itemize}
	\item D\'emonstration de l'\'equation employ\'ee.
	\item Le calcul des racines et la discussion.
	\item Le calcul de $\fii$ et de $h$.
	\item La r\'edaction d'un organigramme et d'un programme en Matlab, accompagn\'es d'une notice d'emploi sur les limites d'emploi, le temps d'\'execution et l'encombrement en machines. 
	\item Un jeu d'essai au moins.
\end{itemize}
\subsubsection{Les M\'ethodes Finies Propos\'ees:}
\begin{figure}
		\centering
			\includegraphics[width=0.50\textwidth]{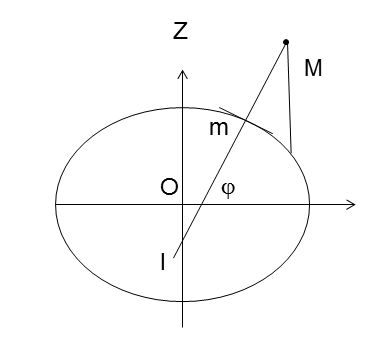}
		\caption{La M\'ethode Finie 1.}
		\label{fig:presentation1s}
	\end{figure}

\begin{figure}
	\centering
	\includegraphics[width=0.7\textwidth]{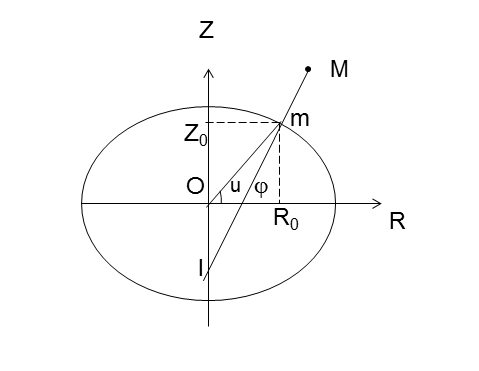}
		\caption{La M\'ethode Finie 2.}
		\label{fig:presentation2s}
	\end{figure}

\begin{enumerate}
	\item On \'ecrit que, simultan\'ement, le point $m$ appartient \`a l'ellipse m\'eridienne et que sa distance \`a $M$ est extr\'emale. C'est un probl\`eme d'extr\'emum li\'e. On se ramenera \`a l'\'equation (\ref{b11}) dans laquelle $x$ est le multiplicateur de Lagrange (Fig.\ref{fig:presentation1s}).
	\item  On \'ecrit que la pente de $mM$ et celle de $mI$ valent $tg\fii=\ds \frac{a^2}{b^2}tgu$ o\`u $u$ est la latitude g\'eocentrique du point $m$. Il vient:
	\be
	\ds \frac{Z-Z_0}{R-R_0}=\ds \frac{a^2}{b^2}.\frac{Z_0}{R_0}
	\ee
	On \'elimine $Z_0$ \`a l'aide de l'\'equation de l'ellipse et on obtient une \'equation du 4\`eme degr\'e en $R_0$ (Fig. \ref{fig:presentation2s}).
	\end{enumerate}
	
	\newpage
\subsection{M\'ethode des D\'eveloppements limit\'es}
$x$ \'etant une inconnue auxiliaire, et $t$ un param\`etre facilement calculable en fonction des donn\'ees, on cherche les $a_i$ tels que:
\be
x=\sum_{i=0}^{i=n}a_it^n
\ee
Les $a_i$ sont des coefficiens fonctions des donn\'ees $X,Y,Z,a$ et $e^2$; ils sont connus jusqu'\`a un certain rang, ou calculables par un programme.

Il peut \^etre difficile de majorer l'erreur de la m\'ethode, suivant les cas. Le d\'eveloppement peut aussi diverger pour certaines valeurs de $t$.
\subsubsection{L'Analyse et la Programmation:}
L'Analyse et la programmation d'une m\'ethode par d\'eveloppement limit\'e peut comporter:
\begin{enumerate}
	\item La d\'emonstration des formules employ\'ees.
	\item Une \'etude analytique de pr\'ecision ou, si impossible, une comparaison avec d'autres m\'ethodes.
	\item La r\'edaction d'un organigramme et d'un programme en Matlab, accompagn\'es d'une notice d'emploi sur les limites d'emploi, le temps d'\'execution et l'encombrement en machines. 
	\item Un jeu d'essai au moins.
\end{enumerate}
\subsubsection{Le D\'eveloppement propos\'e:}
On pose :
\ba
R=(N+h)cos\fii \\
Z=(N(1-e^2)+h)sin\fii \\
x=\ds \frac{R}{Z}tg\fii \\
\nu^2=(1-e^2)Z^2 \\
t=\ds \frac{e^2}{\sqrt{R^2+\nu^2}}\\
\rho=\sqrt{R^2+x^2\nu^2}\\
c=\ds\frac{R^2}{R^2+\nu^2}
\ea
1. Montrer que:
\be
x-1=\ds e^2\frac{ax}{\rho}
\ee
2. Montrer que $x$ v\'erifie l'\'equation:
\be
(x-1)^2(R^2+\nu^2x^2)=e^4a^2x^2
\ee
On pose alors:
\be
x=x=\sum_{i=0}^{i=n}a_it^n
\ee
3. Trouver les $a_i$ pour $0\leq i\leq 2$ par identification des deux membres de l'\'equation en $x$.

On admettra ensuite:
\ba
a_3=\ds \frac{5c^2-3c}{2}\\
a_4=2c-9c^2+8c^3
\ea
			\newpage						
\part{\textsc{Exercices et Probl\`emes}}
\chapter{Exercices et Probl\`emes}
 	\section{	\textsc{Trigonom\'etrie Sph\'erique}}
\textbf{Exercice $n^{\circ}$1:}
 Calculer l'azimut d'une \'etoile de d\'eclinaison $\delta  =+5^{\circ}$ quand sa distance z\'enithale est de $80^{\circ}$ pour un observateur situ\'e \`a la latitude $\varphi   = 56^{\circ}.$ 
\\

\textbf{Exercice $n^{\circ}$2:}
 En appliquant au triangle de position les formules de trigonom\'etrie sph\'erique montrer que l'on peut calculer l'angle horaire $AH_c$ du coucher d'un astre par : $cos AH_c = -tg\varphi.tg\delta $.   
\\

\textbf{Exercice $n^{\circ}$3:}
 Soit un triangle sph\'erique $ABC$. On donne les \'el\'ements suivants:

- $\hat{A}=80.1643\,3\,gr$,

- $\hat{B}=55.7735\,1\,gr$,

- $\hat{C}=64.0626\,1\,gr$,

- $AC= 20.1357\,km$,

- $AB= 22.1435\,km$.

1. Calculer $\alpha=\hat{A}+\hat{B}+\hat{C}$.

2. D\'eterminer $\epsilon$ l'exc\`es sph\'erique de ce triangle.

3. Calculer la fermeture du triangle $ABC$, donn\'ee par:
$$ f=\alpha -200.0000\,0\,gr-\epsilon$$
\\

\textbf{Exercice $n^{\circ}$4:}
 Soit $(\BbS^2)$ une sph\`ere de rayon \'egal \`a 1. Soit un carr\'e sph\'erique $ABCD$ de c\^ot\'e $a$ (arc de grand cercle). On note $\alpha=\hat{A}=\hat{B}=\hat{C}=\hat{D}$. 

1. Montrer que:
$$ cos a=cotg^2\frac{\alpha}{2}$$
2. Donner l'expression de la diagonale $d=l'arc\,AC$.  
\\

\begin{figure}
	\centering
		\includegraphics[width=0.50\textwidth]{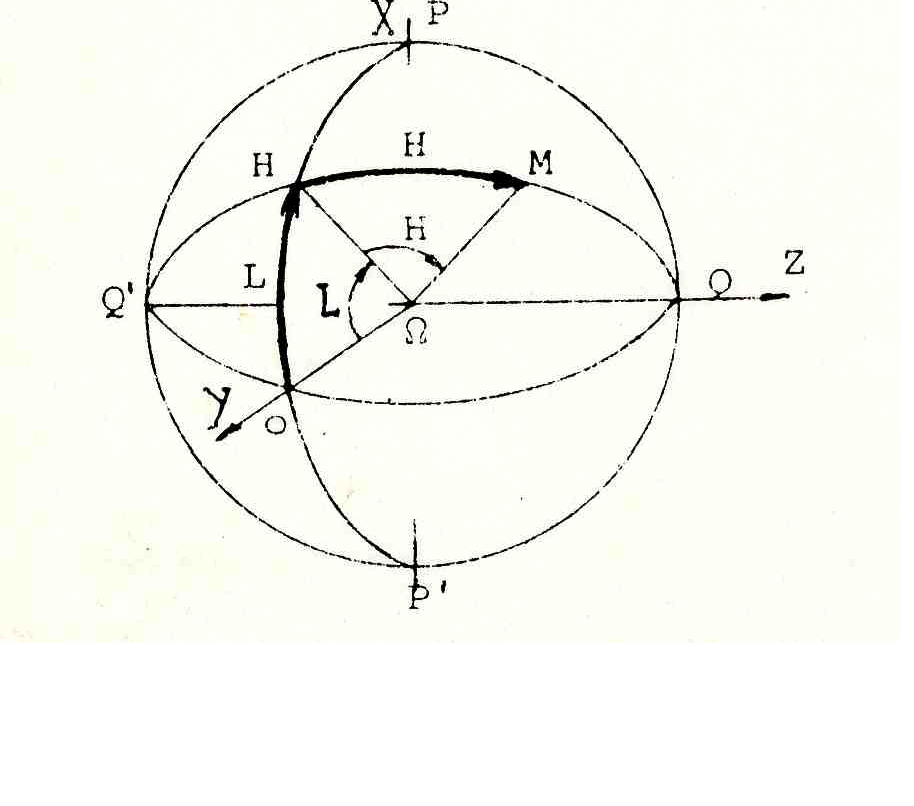}
	\caption{Les coordonn\'ees de Cassini-Soldner}
	\label{fig:soldner}
\end{figure}

\textbf{Probl\`eme $n^{\circ}$1:} 
Soit $(\BbS^2)$ une sph\`ere de rayon \'egal \`a 1 et de centre le point $\Omega$. Un point $M$ de 
$(\BbS^2)$ a pour coordonn\'ees $(\varphi,\lambda)$. On appelle les coordonn\'ees de Cassini-Soldner\index{Coordonn\'ees de Cassini-Soldner}\index{\textbf{Soldner J.G}} \footnote{C\'esar-Fran\c{c}ois Cassini (1714-1784): Astronome et g\'eod\'esien fran\c{c}ais.}\footnote{Dr Johann Georg von Soldner (1776-1833): Math\'ematicien et astronome bavarois.} de $M$ les angles (Fig.\ref{fig:soldner}):\index{\textbf{Cassini C.F.}}

- $L=\widehat{\Omega O,\Omega H}$,

- $H=\widehat{\Omega H,\Omega M}$.

1. D\'eterminer les relations liant $L,H$ \`a $\varphi,\lambda$.

2. Inversement, donner les relations liant  $\varphi,\lambda$ \`a $L,H$.

\textbf{Probl\`eme $n^{\circ}$2:}
 Au lieu $M$ de latitude $\varphi   = 38^{\circ}$  Nord, on observe l'\'etoile polaire  $A$ de d\'eclinaison $\delta   = + 89^{\circ}$ et d'ascension droite $\alpha$  = $+2\,h\, 13\,mn\, 52.90\,s$.
 
1.	Donner sur un graphique, les \'el\'ements du triangle sph\'erique $PAM$ o\`u $P$ est le p\^ole Nord. 

2.	Sachant que l'heure sid\'erale locale $HSL$ est \'egale au moment de l'observation \`a  $6\,h\, 37\, mn\, 19.72\, s$,  calculer l'angle horaire $AH$.  

3.	En appliquant la formule des cotangentes, montrer  que l'azimut $Az$ de l'\'etoile est donn\'e par la formule:
$$tgAz=\frac{sinAH}{cosAHsin\varphi-cos\varphi tg\delta}$$
4.	Calculer alors l'azimut $A$z.

5.	Calculer la distance z\'enithale $z$ de l'\'etoile.
\\
\newpage
\section{	\textsc{Astronomie de Position}}
\textbf{Exercice $n^{\circ}$1:}
 Au lieu de latitude $ \varphi = 36^{\circ} 54'$ Nord, on veut calculer les hauteurs $h_1$ et $h_2$ de l'\'etoile polaire de d\'eclinaison $\delta= + 89^{\circ}$ respectivement \`a son passage sup\'erieur et \`a son passage inf\'erieur au m\'eridien du lieu. D\'eterminer $h_1$ et $h_2$. 
\\ 

\textbf{Probl\`eme $n^{\circ}$1:}
1. En un lieu de latitude $\varphi$ quelles sont les \'etoiles : 

- qui ne se couchent pas ( qui sont toujours visibles),

- qui ne sont jamais visibles. 

Traiter le cas : lieu dans l'h\'emisph\`ere nord. 

2. Quelle est la condition pour qu'une \'etoile culmine au z\'enith ? 

3. Cas particulier du soleil: la d\'eclinaison du soleil varie de $-23^{\circ}27'$ \`a $+23^{\circ}27'$ au cours de l'ann\'ee. On appelle jour le moment pendant lequel le soleil est au-dessus de l'horizon, nuit lorsque le soleil est au-dessous de l'horizon, midi l'instant de la culmination, minuit l'instant du passage inf\'erieur. 

     a) Montrer qu'au moment des \'equinoxes le jour et la nuit sont d'\'egale dur\'ee quel que soit le lieu. 

     b) Montrer qu'\`a l'\'equateur, quelle que soit la date le jour et la nuit sont d'\'egale dur\'ee. 

     c) Au moment du solstice d'hiver quels sont les lieux :

 - o\`u il fait constamment jour,

 - o\`u il fait constamment nuit. 

M\^emes questions au moment du solstice d'\'et\'e. 
   
		d) Quels sont les lieux de la terre o\`u le soleil culmine au z\'enith au moment du solstice d'hiver. M\^eme question au moment du solstice d'\'et\'e. 
    
		e) Quels sont les lieux de la terre o\`u au moins une fois dans l'ann\'ee le soleil culmine au z\'enith. 
\\

\textbf{Probl\`eme $n^{\circ}$2:}
Une station astronomique est situ\'ee en un lieu de coordonn\'ees g\'eographiques : $\varphi =\, + 45^{\circ}\,00';\,\, \lambda =\, + 7\, h\, 20\, mn$.

En ce lieu, on observe une \'etoile $A$ de coordonn\'ees \'equatoriales: 
 
$\alpha = \,+11\, h\, 13\, mn;\,\,\delta =\, 30^{\circ}\,00'.$

L'observation se fait le jour de l'\'equinoxe de printemps le 21 mars \`a $0$ heure $TU.$ L'heure sid\'erale de Greenwich est $11\, h\, 52\, mn$. 

1. Calculer l'heure sid\'erale locale du lever et du coucher de l'\'etoile $A$ au lieu consid\'er\'e.
 
2. En d\'eduire l'heure $TU$ du lever et du coucher de l'\'etoile au lieu consid\'er\'e. 

Remarque: on choisira le coucher qui a lieu apr\`es le lever. 
\\

\textbf{Probl\`eme $n^{\circ}$3:}
 En un lieu de latitude $43^{\circ},521$ et de longitude $+ 0\,h\,20\,mn\,57\,s$, on cherche \`a pointer la galaxie d'Androm\`ede de coordonn\'ees \'equatoriales $\alpha = 0\,h\,40\,mn,\, \delta = 41^{\circ}\,00'$ le 31 juillet 1992 \`a $21\, h\, TU.$ 

On donne l'heure sid\'erale de Greenwich \`a $0\,h\,TU$ le 31/07/1992: $HSG_{0hTU} = 20\,h\,35\,mn\,28\,s.$ 

1. Calculer l'heure sid\'erale locale \`a $21\, h\, TU.$

2. En d\'eduire l'angle horaire de la galaxie.

3. Calculer la distance z\'enithale de la galaxie \`a $21\, h\, TU.$

4. Calculer son azimut \`a cette m\^eme heure.
\\

\textbf{Probl\`eme $n^{\circ}$4:}
En un lieu de l'h\'emisph\`ere Nord de latitude $\varphi$, on mesure la longueur de l'ombre port\'ee $HC$, \`a midi vrai (passage du soleil au m\'eridien), par une tige verticale $HA$ dont l'extr\'emit\'e $H$ est sur le sol suppos\'e horizontal. 

1. Donner l'expression $HC$ en fonction de $HA$ et de la distance z\'enithale $Dz$ du soleil. 

2. Donner l'expression de $HC$ en fonction de $HA$ et de $\varphi$: 

- aux \'equinoxes,

- aux solstices.
	
3. Quelle doit \^etre la d\'eclinaison du soleil et en quels lieux, pourque l'on ait $HC = HA$? 

4. En un lieu de latitude $\varphi = 47^{\circ}$ en quelles saisons peut on avoir $HC = HA$.
 
5. Si on d\'eplace $HB$ le long d'un m\'eridien, en restant dans l'h\'emisph\`ere Nord, existe-t- il au cours de l'ann\'ee des lieux o\`u $HC= 0$, ou $HC$ devient infiniment grand.

\newpage
\section{	\textsc{Courbes et Surfaces}}
\textbf{Exercice $n^{\circ}$1:}
Soit l'h\'elice circulaire $\Gamma$ param\'etr\'ee par:
$$\left\{ \begin{array}{l}
x=acost\\
y=asint\\
z=bt
\end{array}\right. $$
o\`u $a,b$ deux constantes positives.

1. Exprimer les composantes des vecteurs $T,N,B$ du rep\`ere de Fr\'enet.

2. Montrer que la courbure vaut $\ds \frac{a}{a^2+b^2}$.

3. Montrer que la torsion vaut $\ds \frac{b}{a^2+b^2}$.
\\

\textbf{Probl\`eme $n^{\circ}$1:}
 Soit la courbe $(C)$ d\'efinie par les formules:
$$M \left\{\begin{array}{l}
x=at^2 \\
y= at^3 \\
z=\ds \frac{9}{16}at^4\quad \mbox{avec}\quad a>0
\end{array}\right. $$
1. Calculer l'abscisse curviligne $s$ d'un point $M$ quelconque de cette courbe lorsqu'on prend pour origine des arcs l'origine des coordonn\'ees et qu'on prend pour sens des arcs croissants celui des $y$ croissants.

2. D\'eterminer au point $M$ les vecteurs unitaires du tri\`edre de Fren\^et.

3. Calculer le rayon de courbure et les coordonn\'ees du centre de courbure.

4. Evaluer la torsion en $M$.
\\

\textbf{Probl\`eme $n^{\circ}$2:}
Soit $(\Gamma)$ la surface param\'etr\'ee par $(u,v)$ dans $\BbR^2$ telle que:
\[M(u,v) \left\{ \begin{array}{l}
	     X =u(1-u^2)cosv \\
	    Y = u(1-u^2)sinv\\
	    Z = 1-u^2
	    \end{array} \right.
\]
1. Calculer l'expression de $ds^2$.

2. Montrer que l'\'equation cart\'esienne de $(\Gamma)$ est: 
$$ x^2+y^2=(1-z)z^2 $$
\\

\textbf{Exercice $n^{\circ}$2:}
 Soit la surface d'Enneper:
	\[M(u,v) \left\{ \begin{array}{l}
	     X =\ds u-\frac{u^3}{3}+uv^2 \\
	    Y = v-\ds \frac{v^3}{3}+vu^2 \\
	    Z = u^2-v^2
	    \end{array} \right.
\]
1. Montrer que:    $$ ds^2 = (1+u^2+v^2)^2.(du^2+dv^2)$$

2. Calculer un vecteur unitaire normal \`a la surface.

3. Montrer que la surface d'Enneper est de courbure moyenne nulle en chaque point.
\\

\textbf{Exercice $n^{\circ} 3$ :}
 On suppose que la m\'etrique d'une surface donn\'ee est:
$$ds^2=A^2du^2+B^2dv^2, \quad A=A(u,v),\quad B=B(u,v)$$
1. Montrer alors que l'expression de la courbure totale est:
$$ K=-\frac{1}{AB}\left[\left(\frac{A'_v}{B}\right)'_v+\left(\frac{B'_u}{A}\right)'_u\right]$$
$'$ d\'esigne la d\'erivation partielle. 
\\

\textbf{Probl\`eme $n^{\circ}$3:}
 Soit l'ellipse $(E)$ d\'efinie par les \'equations param\'etriques:
$$ M \left\{\begin{array}{l}
x=acos u \\
y= bsinu \\
\mbox{avec}\quad a>b>0
\end{array}\right. $$
On  pose:
$$e^2=\frac{a^2-b^2}{a^2};\quad e'^2= \frac{a^2-b^2}{b^2}$$
1. Calculer la position sur l'axe des abscisses des deux points $F$ et $F'$ appel\'es foyers tels que $MF+MF'=2a$.

2. Montrer que le produit des distances des foyers \`a la tangente \`a l'ellipse en M est ind\'ependant de $u$.

3. Donner l'expression de $ds$.

4. D\'eterminer les expressions des vecteurs unitaires $\textbf{\textit{T}}$ et $\textbf{\textit{N}}$ et en d\'eduire le rayon de coubure de l'ellipse.  

5. Montrer qu'il passe par $M$ deux cercles tangents en ce point \`a la courbe et centr\'es sur $Ox,Oy$ respectivement (appel\'es cercles surosculateurs).

6. Que deviennent ces cercles lorsque $M$ est un sommet de l'ellipse.
\\

\textbf{Probl\`eme $n^{\circ}$4:} 
On d\'efinit une surface $(S)$ par les \'equations:
	\[M(u,v) \left\{ \begin{array}{l}
	     X = u^2 + v \\
	    Y = u + v^2 \\
	    Z = uv
	    \end{array} \right.
\]
1.	Calculer les composantes des vecteurs $\textit{\textbf{OM}}'_u$ et $\textit{\textbf{OM}}'_v$.

2.	Calculer les coefficients $E,F,G$ de la premi\`ere forme fondamentale de la surface $(S)$.

3.	En d\'eduire l'expression de $ds^2$.

4.	Les coordonn\'ees $(u,v )$ sont-elles orthogonales? sym\'etriques?

5.	Calculer un vecteur normal de $(S)$. 
\\

\textbf{Probl\`eme $n^{\circ}$5:} 
On d\'efinit une surface $(\Sigma)$ par les \'equations:
 	\[M(u,v) \left\{ \begin{array}{l}
	           X =  a.cosu.cosv \\
	         Y = a.cosu.sinv \\
	          Z = b.sinu
	              \end{array} \right.
\] 
avec $a,b$ deux constantes positives. 

1.	Calculer les composantes des vecteurs $\textit{\textbf{OM}}'_u$ et $\textit{\textbf{OM}}'_v$.

2.	Calculer les coefficients $E,F,G$ de la premi\`ere forme fondamentale de la surface $(\Sigma)$.

3.	En d\'eduire l'expression de $ds^2$.

4.	Les coordonn\'ees $( u,v )$ sont-elles orthogonales? sym\'etriques?

5.	Calculer un vecteur unitaire normal $\textit{\textbf{n}}$ de $(\Sigma)$.

6.	Calculer les vecteurs :
$$   \textit{\textbf{OM}}''_{uu},\quad \textit{\textbf{OM}}''_{uv},\quad \textit{\textbf{OM}}''_{vv}$$
On pose:
$$L = \textit{\textbf{n}}.\textit{\textbf{OM}}''_{uu},\quad  M = \textit{\textbf{n}}.\textit{\textbf{OM}}''_{uv},\quad N = \textbf{\textit{n}}.\textit{\textbf{OM}}''_{vv} $$
7. Calculer les coefficients $L,M$ et $N$.
\\

\textbf{Probl\`eme $n^{\circ}$6:}
 On consid\`ere la surface $(\Gamma)$ d\'efinie par les \'equations:
 	\[M(u,v)\left\{ \begin{array}{lll}
	           X =  sinu.cosv \\
	         Y = sinu.sinv \\
	          Z = cosu+Logtg\ds\frac{u}{2} +\psi(v)
	              \end{array} \right.
\] 
avec $\psi(v)$ est une fonction d\'efinie de classe $C^1$ de $v$. 

1. Donner le domaine de d\'efinition de la surface  $(\Gamma)$.

2.	Montrer que les courbes coordonn\'ees $v=constante$ constituent une famille de courbes planes de $(\Gamma)$ et que leur plan coupe $(\Gamma)$ sous un angle constant.

3.  Calculer les composantes des vecteurs $\textit{\textbf{OM}}'_u$ et $\textit{\textbf{OM}}'_v$.

4.	Calculer les coefficients $E,F,G$ de la premi\`ere forme fondamentale de la surface $(\Gamma)$.

5.	En d\'eduire l'expression de $ds^2$.

6.	Les coordonn\'ees $( u,v )$   sont-elles orthogonales? sym\'etriques?

7.	On suppose pour la suite que $\psi(v)=0$, calculer un vecteur unitaire normal $\textit{\textbf{n}}$ de $\Gamma$.

8.	Calculer les vecteurs :
$$   \textit{\textbf{OM}}''_{uu},\quad \textit{\textbf{OM}}''_{uv},\quad \textit{\textbf{OM}}''_{vv}$$
On pose:
$$L = \textit{\textbf{n}}.\textit{\textbf{OM}}''_{uu},\quad  M = \textit{\textbf{n}}.\textit{\textbf{OM}}''_{uv},\quad N = \textbf{\textit{n}}.\textit{\textbf{OM}}''_{vv} $$
9. Calculer les coefficients $L,M$ et $N$.

10. En d\'eduire l'expression des courbures moyenne et totale.
\\

\textbf{Probl\`eme $n^{\circ}$7:}
 Soit la surface $(\Gamma)$ d\'efinie param\'etriquement par:
 	\[M(u,v)\left\{ \begin{array}{lll}
	           X =  thu.cosv \\
	         Y = thu.sinv \\
	          Z = \ds \frac{1}{chu}+Logth\ds\frac{u}{2} 
	              \end{array} \right.
\] 
avec $chu$ et $thu$ sont respectivement le cosinus et la tangente hyperboliques d\'efinies par:
$$ chu=\frac{e^u+e^{-u}}{2},\quad thu=\frac{e^u+e^{-u}}{e^u-e^{-u}}$$ 
1.  Donner le domaine de d\'efinition de la surface  $(\Gamma)$.

2.  Calculer les composantes des vecteurs $\textit{\textbf{OM}}'_u$ et $\textit{\textbf{OM}}'_v$.

3.	Calculer les coefficients $E,F,G$ de la premi\`ere forme fondamentale de la surface $(\Gamma)$.

4.	En d\'eduire l'expression de $ds^2$.

5.	Les coordonn\'ees $( u,v )$   sont-elles orthogonales? sym\'etriques?

6.	Calculer un vecteur unitaire normal $\textit{\textbf{n}}$ de $(\Gamma)$.

7.	Calculer les vecteurs :
$$   \textit{\textbf{OM}}''_{uu},\quad \textit{\textbf{OM}}''_{uv},\quad \textit{\textbf{OM}}''_{vv}$$

On pose:
$$L = \textit{\textbf{n}}.\textit{\textbf{OM}}''_{uu},\quad  M = \textit{\textbf{n}}.\textit{\textbf{OM}}''_{uv},\quad N = \textbf{\textit{n}}.\textit{\textbf{OM}}''_{vv} $$
8. Calculer les coefficients $L,M$ et $N$.

9. D\'eterminer les coubures moyenne et totale.
\\

\textbf{Probl\`eme $n^{\circ}$8:}
 Montrer que les courbures totale $K$ et moyenne $H$ en un point $M(x,y,z)$ d'une surface param\'etr\'ee par $z=f(x,y)$, o\`u $f$ est une fonction lisse, sont donn\'ees par:
$$K=\ds \frac{f''_{xx}f''_{yy}-f''^2_{xy}}{(1+f'^2_x+f'^2_y)^2}$$
 et:
$$ \ds H=\frac{(1+f'^2_x)f''_{xx}-2f'_xf'_yf''_{xy}+(1+f'^2_x)f''_{yy}}{(1+f'^2_x+f'^2_y)^{\frac{3}{2}}} $$
\\

\textbf{Probl\`eme $n^{\circ}$8:}
 Soit $(\Sigma)$ une surface de $\BbR^3$ param\'etr\'ee par $OM(u,v)$ telle que sa premi\`ere forme fondamentale s'\'ecrit: $ds^2=Edu^2+2Fdudv+Gdv^2$
		
1. Montrer que les conditions suivantes sont \'equivalentes:

\quad i) - $\ds \frac{\partial E}{\partial v}=\frac{\partial G}{\partial u}=0$,
	
\quad	ii) - Le vecteur $\ds \frac{\partial^2 OM }{\partial u\partial v}$ est parall\`ele au vecteur normal $N$ \`a la surface,

\quad iii) - Les c\^ot\'es oppos\'es de tout quadrilat\`ere curviligne form\'es par les courbes coordonn\'ees $(u,v)$ ont m\^eme longueurs.

2. Quand ces conditions sont satisfaites, on dit que les courbes coordonn\'ees de $(\Sigma)$ forment un r\'eseau de \textit{Tchebychev}.\footnote{Pafnouti Tchebychev (1821 - 1894 ): Math\'ematicien russe.} Montrer que dans ce cas, on peut param\'etrer la surface par $(\tilde{u},\tilde{v})$ telle que $ds^2$ s'\'ecrit:\index{\textbf{Tchebychev P.}}
$$ ds^2=d\tilde{u}^2+2cos\theta d\tilde{u}d\tilde{v}+d\tilde{v}^2$$
o\`u $\theta$ est une fonction de $(\tilde{u},\tilde{v})$. Montrer que $\theta$ est l'angle entre les courbes coordonn\'ees $\tilde{u},\tilde{v}$.

3. Montrer que l'expression de la courbure totale est donn\'ee par:
$$ K= \ds \frac{1}{sin\theta}.\ds \frac{\partial^2\theta}{\partial \tilde{u} \partial \tilde{v}}$$
4. On pose :
$$\begin{array}{l}
\hat{u}=\tilde{u}+\tilde{v} \\
\hat{v}=\tilde{u}-\tilde{v}
\end{array}$$
Montrer que $ds^2$ s'\'ecrit avec les nouvelles variables $(\hat{u},\hat{v})$:
$$ds^2=cos^2\omega d\hat{u}^2+sin^2\omega d\hat{v}^2$$ avec $\omega=\theta/2$. (\textit{A.N. Pressley}, 2010)\index{\textbf{Pressley A.N.}}

\newpage
\section{	\textsc{La G\'eom\'etrie de l'Ellipse et de l'Ellipsoide}}
\textbf{Exercice $n^{\circ}$1:} 
A partir de la d\'efinition g\'eom\'etrique de l'ellipse donn\'ee par: $$ MF + MF' = constante=2a $$ retrouver l'expression de l'\'equation cart\'esienne de l'ellipse.
\\

\textbf{Exercice $n^{\circ}$2:}
 Montrer la formule tr\`es utilis\'ee en g\'eod\'esie:
$$ \frac{d(Ncos\varphi)}{d\varphi}=-\rho sin\varphi $$
 avec $N$ et $\rho$ les deux rayons de courbures principaux de l'ellipsoide de r\'evolution donn\'es respectivement par :
$$ N=\frac{a}{\sqrt{1-e^2sin^2\fii}}$$
et: $$ \rho=N=\frac{a(1-e^2)}{(1-e^2sin^2\fii)\sqrt{1-e^2sin^2\fii}}$$ 
\\

 \textbf{Probl\`eme $n^{\circ}$1:} 
A partir des \'equations de l'ellipsoide de r\'evolution:
$$M=\left\{\begin{array}{lll}
	X=Ncos\varphi cos\lambda  \\
	Y = Ncos\varphi sin\lambda  \\
	Z=N(1-e^2)sin\varphi   \end{array}\right.
	$$
1. Calculer les vecteurs:
$$ \displaystyle \frac{\partial \textbf{\textit{M}}}{\partial \lambda},\frac{\partial \textbf{\textit{M}}}{\partial \varphi}$$
2. Calculer les coefficients:
$$ E=\displaystyle \frac{\partial \textbf{\textit{M}}}{\partial \lambda}.\frac{\partial \textbf{\textit{M}}}{\partial \lambda}, \quad
	F=\displaystyle \frac{\partial \textbf{\textit{M}}}{\partial \lambda}.\frac{\partial \textbf{\textit{M}}}{\partial \varphi}, \quad
		G=\displaystyle \frac{\partial \textbf{\textit{M}}}{\partial \varphi}.\frac{\partial \textbf{\textit{M}}}{\partial \varphi}
   	$$
	D\'emontrer que l'expression de la premi\`ere forme fondamentale s'\'ecrit:
$$ ds^2=\rho^2d\varphi^2+N^2cos^2\varphi d\lambda^2 $$
3. Calculer le vecteur normal $\textbf{\textit{n}}$ :
\begin{displaymath}
\textbf{\textit{n}}=\displaystyle \frac{\partial \textbf{\textit{M}}}{\partial \lambda} \wedge \frac{\partial \textbf{\textit{M}}}{\partial \varphi}\frac{1}{\left\|\displaystyle \frac{\partial \textbf{\textit{M}}}{\partial \lambda} \wedge \frac{\partial \textbf{\textit{M}}}{\partial \varphi}\right\|}
\end{displaymath}
4. Calculer les vecteurs:
$$ 	\displaystyle \frac{\partial^2 \textbf{\textit{M}}}{\partial \lambda^2},  \quad	\displaystyle \frac{\partial^2 \textbf{\textit{M}}}{\partial \lambda \partial \varphi}, \quad  \displaystyle \frac{\partial^2 \textbf{\textit{M}}}{\partial \varphi^2} $$
5. D\'eterminer les coefficients:
$$  	L=\displaystyle n.\frac{\partial^2 \textbf{\textit{M}}}{\partial \lambda^2}, \quad 
	M=\displaystyle n.\frac{\partial^2 \textbf{\textit{M}}}{\partial \lambda \partial \varphi}, \quad
		N=\displaystyle n.\frac{\partial^2 \textbf{\textit{M}}}{\partial ^2\varphi} 	$$
6. Ecrire la deuxi\`eme forme fondamentale $\Phi(\lambda,\varphi)$.

7. En appliquant la formule du cours, Montrer que :
\begin{displaymath}
  N(\varphi)=\frac{a}{\sqrt{1-e^2sin^2\varphi}}
\end{displaymath}
est le rayon de courbure de la section normale au point $M$ perpendiculaire au plan de la m\'eridienne de l'ellipsoide de r\'evolution.

8. En posant:
\begin{displaymath}
  d\m L= \frac{\rho d\varphi}{Ncos\varphi}
\end{displaymath}
En d\'eduire que $ds^2$ s'\'ecrit:
\begin{displaymath}
  ds^2=N^2cos^2\varphi (d\m L^2+d\lambda^2)
  \end{displaymath}
9. Montrer que $\mathcal L$ est donn\'ee par:
\begin{displaymath}
\m L(\varphi)=Log\left(tg(\frac{\pi}{4}+\frac{\varphi}{2})\right)-\frac{e}{2}Log\left(\frac{1+esin\varphi}{1-esin\varphi}\right)
\end{displaymath}
\\

\textbf{Probl\`eme $n^{\circ}$2:}
 Sur l'ellipsoide, on note $\varphi$  la latitude g\'eod\'esique et  $\psi$   la latitude r\'eduite.

1. Calculer $\rho$  le  rayon de  courbure de l'ellipse m\'eridienne en fonction de $\psi$.

2. Exprimer l'aplatissement de l'ellipsoide en fonction des valeurs de $\rho$  au p\^ole et \`a l'\'equateur. 

3. On mesure la longueur d'un arc de m\'eridien d'un degr\'e \`a la fois au p\^ole et \`a l'\'equateur. On trouve respectivement $111\,695\, m$ et $110\,573\, m$. En d\'eduire l'aplatissement.
\\

\textbf{Probl\`eme $n^{\circ}$3:}
 On donne les coordonn\'ees tridimensionnelles suivantes d'un point $M$:
	\[	M=(X,Y,Z)=(4\,300\,244.860\, m,1\,062\,094.681\, m,4\,574\,775.629\, m)
\]
Les param\`etres de l'ellipsoide de r\'ef\'erence sont $a  = 6\,378\,137.00\, m,\quad  e^2 = 0.006\,694\,38$.
 
1.	Calculer le demi-petit axe $b$.

2.	Calculer l'aplatissement. 

3.	Calculer les coordonn\'ees g\'eod\'esiques $(\varphi ,\lambda , he)$ du point $M$. $\varphi$   et $\lambda$    seront calcul\'ees en grades avec  cinq chiffres apr\`es la virgule.
\\

\textbf{Probl\`eme $n^{\circ}$4:} 
Soit $\m E(a,e)$ un ellipsoide de r\'evolution o\`u $a,e$ sont respectivement le demi-grand axe et la premi\`ere excentricit\'e. $(g)$ une g\'eod\'esique partant d'un point $E(\varphi=0,\lambda_E)$ sur l'\'equateur et d'azimut $Az_E$. A cette g\'eod\'esique, on lui fait correspondre une g\'eod\'esique $(g')$ sur la sph\`ere $\m S^2$ dite de Jacobi\footnote{Carl Gustav Jacob Jacobi (1804-1851): Math\'ematicien allemand.} \index{\textbf{Jacobi C.G.J.}} de rayon $a$, ayant le m\^eme azimut $Az_E$ au point $E'(\varphi'=0,\lambda_E)$. De m\^eme au point $M(\varphi,\lambda)$ de la g\'eod\'esique $(g)$ de l'ellipsoide, on lui fait correspondre le point $M'(\varphi',\lambda')$ de $(g')$ de $\m S^2$ tel qu'il y a conservation des azimuts. 
\begin{figure}[h]
	\centering
		\includegraphics[width=0.70\textwidth]{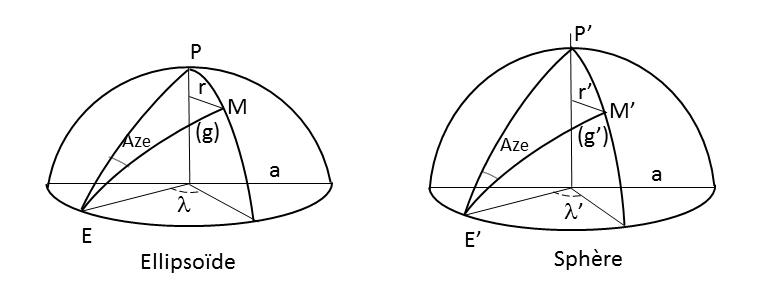}
	\caption{La Correspondance de la sph\`ere de Jacobi}
	\label{fig:sjacobi}
\end{figure}

1. Ecrire l'\'equation de Clairaut pour la g\'eod\'esique $(g)$. 

2. On note $r'$ le rayon du parall\`ele passant par $M'$ de la g\'eod\'esique $(g')$. Ecrire de m\^eme l'\'equation de Clairaut pour la g\'eod\'esique $(g')$. 

3. Montrer que $\varphi$ et $\varphi'$ v\'erifient:
$$ Ncos\varphi=acos\varphi'$$
et en d\'eduire que $\varphi'$ est la latitude param\'etrique de $M$.

4. Ecrire les expressions de $tgAz_g$ et $tgAz_{g'}$ respectivement sur $(g)$ et $(g')$.

5. Montrer que:
$$  d\lambda=\frac{\rho d\varphi}{ad\varphi'}d\lambda'$$ 
En d\'eduire que : 
$$ 	d\lambda=\sqrt{1-e^2cos^2\varphi'}d\lambda' $$
6. En int\'egrant l'\'equation pr\'ec\'edente, montrer qu'on obtient:
$$	\lambda - \lambda_E=\int_{\lambda_E}^{\lambda'+\lambda_E}\sqrt{1-e^2cos^2\varphi'}d\lambda' $$
avec $\lambda>\lambda_E$ et $\lambda'$ est compt\'ee \`a partir de $\lambda_E$.

7. En \'ecrivant $ \sqrt{1-e^2cos^2\varphi'}=1-\frac{e^2}{2}cos^2\varphi'+o(e^4)$ o\`u $o(e^4)$ est un infiniment petit d'ordre 4 en $e$ dont on n\'eglige, \'ecrire l'int\'egrale pr\'ec\'edente entre $\lambda_E$ et $\lambda_E+\lambda$.

8. Comme $(g')$ est une g\'eod\'esique de la sph\`ere, on d\'emontre que:
$$ 	cos^2\varphi'd\lambda'=\frac{sinAz_E}{a}ds'  $$
o\`u $ds'$ est l'\'el\'ement diff\'erentiel de l'abscisse curviligne sur la g\'eod\'esique (un grand cercle). Alors en posant $s'=0$ au point $E'$, montrer que l'\'equation pr\'ec\'edente s'\'ecrit sous la forme:
$$ 		\lambda =\lambda_E+\lambda'-\frac{e^2sinAz_e}{2a}\int_0^{s'}ds' $$ 
9. On suppose que la g\'eod\'esique $(g')$ coupe une premi\`ere fois le plan de l'\'equateur en un point $F'$, montrer qu'on obtient:
\ba
	\lambda'_F=\pi \label{27}\nonumber \\
	s'=\pi a \label{28}  \nonumber \\
	\lambda_F=\lambda_E+\pi-\frac{e^2\pi sinAz_E}{2}  \nonumber
\ea
10. La g\'eod\'esique $(g')$ partant de $F'$ a pour azimut $\pi-Az_E$, elle coupe une deuxi\`eme fois l'\'equateur au point $E'$, mais la g\'eod\'esique $(g)$ sur l'ellipsoide coupe une deuxi\`eme fois le plan de l'\'equateur au point correspondant \`a $H$ dont la longitude est $\lambda_H$. Montrer que $\lambda_H$ est donn\'ee par:
$$	\lambda_H=\lambda_E+2\pi-\frac{e^2\pi sinAz_E}{2}-\frac{e^2\pi sin(\pi-Az_E)}{2}=\lambda_E+2\pi-e^2\pi sinAz_E $$
Quelle conclusion a-t-on sur les lignes g\'eod\'esiques de l'ellipsoide de r\'evolution.
  \\

\textbf{Probl\`eme $n^{\circ}$5:} 
Un point $M$ de la surface d'une sph\`ere $(S)$ de rayon $R$, a pour coordonn\'ees $(X,Y,Z)$ dans un rep\`ere orthonorm\'e: 
$$	  M=(X,Y,Z)=(Rcos\varphi .cos\lambda ,Rcos\varphi .sin \lambda ,Rsin\varphi)$$
1. Montrer qu'un vecteur normal unitaire $n$ \`a $(S)$ en $M$ est:
$$ n=	 (cos\varphi .cos\lambda,	 cos\varphi .sin \lambda,	 sin\varphi)^T$$
2. Soit $(C)$ le grand cercle passant par le point $A(R,0,0)$ et d'azimut $Az_E$. Le point $M$ peut \^etre d\'ecrit par son abscisse curviligne $s$ mesurant l'arc $AM$. On note par $\omega$ repr\'esente l'angle au centre de l'arc $AM$. Utilisant la trigonom\'etrie sph\'erique, montrer que:
	$$cos\varphi.sin\lambda= sin\omega.sinAz_E $$
3. En utilisant la formule fondamentale de la trigonom\'etrie sph\'erique dans le triangle APM, montrer qu'on a les deux  relations :
$$\begin{array}{l}
	      cos\omega =  cos\varphi.cos\lambda \\
	                  sin\varphi = sin\omega.cosAz_E
										\end{array}$$
4. En d\'eduire que les coordonn\'ees de $M$ s'\'ecrivent en fonction de $s$ comme suit:
$$ M\left\{
\begin{array}{l}
	 X= R.cos(s/R)\\
	 Y=RsinAz_Esin(s/R)\\
	 Z=RcosAz_Esin(s/R)
\end{array}\right.$$
5. Calculer les vecteurs $T$ et $N$ du rep\`ere de Fren\^et. En d\'eduire les composantes de $N$ en fonction de $\omega$.

6. Montrer que les vecteurs $N$ et $n$ sont parall\`eles.

7. Justifier que les g\'eod\'esiques de la sph\`ere sont les grands cercles. 
\\

\textbf{Probl\`eme $n^{\circ}$6:} 
Soit le tore $\BbT$ d\'efini par les \'equations suivantes:
$$M(\fii,\lm)=\left\{\begin{array}{l}
x=(a+Rcos\fii)cos\lm \\
y=(a+Rcos\fii)sin\lm \\
z=Rsin\fii 
\end{array}\right.
$$
o\`u $a,R$ deux constantes positives avec $a>R$, $(\fii,\lm)\in \,[0,2\pi]\times [0,2\pi]$. 

1. Calculer la premi\`ere forme fondamentale $ds^2$.

2. Avec les notations usuelles, on pose: 
$$\frac{\partial E}{\partial \fii}=E'_\fii,\quad	\frac{\partial E}{\partial \lm}  =E'_\lm,\quad 	\frac{\partial F}{\partial \fii} =F'_\fii$$
 $$\quad \frac{\partial F}{\partial \lm} =F'_\lm,\quad \frac{\partial G}{\partial \fii}=G'_\fii,\quad	\frac{\partial G}{\partial \lm}= G'_\lm $$
Utilisant les \'equations des g\'eod\'esiques du cours, montrer que les \'equations des g\'eod\'esiques du tore sont:
\[-2Rsin\fii(a+Rcos\fii)	\frac{d\fii}{ds} \frac{d\lm}{ds}+(a+Rcos\fii)^2\frac{d^2 \lm}{ds^2}=0 \]
\[Rsin\fii(a+Rcos\fii)\left(\frac{d\lm}{ds}\right)^2+R^2\frac{d^2 \fii}{ds^2}=0 \]	
3.  Montrer que la premi\`ere \'equation ci-dessus donne: $$ \ds (a+Rcos\fii)^2\frac{d\lm}{ds}=C=cte$$
Montrer qu'on retrouve l'\'equation de Clairaut avec $C=(a+R)sinAze$ o\`u $Aze$ est l'azimut de d\'epart au point $M_0(\fii=0,\lm_0)$.

4.  On suppose au point $M_0$, la g\'eod\'esique a pour azimut $Aze$ tel que:
 $$ 0 < Aze< \ds \frac{\pi}{2}$$
 Montrer que la deuxi\`eme \'equation des g\'eod\'esiques s'\'ecrit en utilisant le r\'esultat pr\'ec\'edent:
	\[  \frac{d^2 \fii}{ds^2}=-\frac{C^2}{R}\frac{sin\fii}{(a+Rcos\fii)^3}
\]
5.  Montrer qu'on arrive \`a: 	$$  \ds \left(\frac{d\fii}{ds}\right)^2=l-\frac{C^2}{R^2(a+Rcos\fii)^2}\geq 0$$
 o\`u $l$ est une constante d'int\'egration.

\newpage
\section{	\textsc{Les Syst\`emes G\'eod\'esiques}}

\textbf{Exercice $n^{\circ}$1:}
 Donner l'expression des composantes du gradient en coordonn\'ees cylindriques.
\\

\textbf{Exercice $n^{\circ}$2:} 
On donne l'expression scalaire d'une fonction $V(x,y,z)$  par : 
 	\[V(x,y,z)=\frac{ax^2+y^2}{z^2}+\frac{1}{2}\omega^2(x^2+y^2)
\]
1. Calculer les composantes du  vecteur $\textit{\textbf{grad}}V$ dans un domaine de $\BbR^3$ o\`u $z\neq 0$. 
\\

\textbf{Probl\`eme $n^{\circ}$1:} 
Soit un point $ A (\varphi ,\lambda)$ sur un ellipsoide de r\'evolution associ\'e \`a un r\'ef\'erentiel g\'eocentrique donn\'e $ \m R$. On consid\`ere le rep\`ere orthonorm\'e local en A $(e_{\lambda},e_{\varphi},e_n)$ d\'efini dans la base orthonorm\'ee $(i,j,k)$ de $ \m R$ o\`u $e_{\lambda}$ est tangent au parall\`ele passant par $A$ et dirig\'e vers l'Est, $e_{\varphi}$ tangent \`a la m\'eridienne, dirig\'e vers le nord et $e_n$ port\'e par la normale \`a l'ellipsoide dirig\'e vers le z\'enith.

1. Exprimer les vecteurs de la base $(e_{\lambda},e_{\varphi},e_n)$ dans la base $(i,j,k)$ de $ \m R$.

2. Exprimer les vecteurs $i,j$ et $k$ dans la base $(e_{\lambda},e_{\varphi},e_n)$.

3. Calculer $de_{\lambda},de_{\varphi}$ et $de_n$ dans la base $(i,j,k)$.

4. En adoptant une \'ecriture matricielle, montrer que :
$$	\begin{pmatrix}
	de_{\lambda} \\
de_{\varphi}  \\
	de_n 
\end{pmatrix}=\begin{pmatrix}
0 &	sin\varphi d\lambda & -cos\varphi d\lambda \\
	-sin\varphi d\lambda & 0 & -d\varphi  \\
	cos\varphi d\lambda & d\varphi & 0
\end{pmatrix}\begin{pmatrix}
	e_{\lambda} \\
	e_{\varphi} \\
	e_n 
\end{pmatrix}$$
\\

\textbf{Probl\`eme $n^{\circ}$2:}
  On d\'efinit dans $\BbR^3$ un point $M$ par ses coordonn\'ees ellipsoidiques de Jacobi $(\phi,\lambda,u)$ comme suit:
$$ M \left\{ \begin{array}{l}
x=\sqrt{u^2+\epsilon^2}.cos\phi cos\lambda \\
y=\sqrt{u^2+\epsilon^2}.cos\phi sin\lambda \\
z=u.sin\phi 
\end{array} \right. $$
avec: $\epsilon^2=\sqrt{a^2-b^2},\,\phi \in [-\pi/2,\pi/2 ],\, \lambda \in [0,2\pi]$ et $u\in ]0,+\infty[$, $a,b$ deux constantes r\'eelles telles que $a>b>0$. 

1. Montrer que le point $M$ appartient \`a un ellipsoide de r\'evolution en pr\'ecisant ses demi-axes.

2. Calculer $ds^2$ et montrer qu'il s'\'ecrit sous la forme:
$$ ds^2=(d\phi,d\lambda,du).G.\begin{pmatrix}
d\phi \\
d\lambda \\
du 
\end{pmatrix}$$
avec $G$ donn\'ee par :
$$ G=(g_{ij})=\begin{pmatrix}
u^2+\epsilon^2sin^2\phi & 0 & 0 \\ 
0 & (u^2+\epsilon^2)cos^2\phi & 0  \\
0 & 0 & \ds \frac{u^2+\epsilon^2sin^2\phi}{u^2+\epsilon^2} 
\end{pmatrix}$$

3. Sachant que l'expression du laplacien d'une fonction scalaire $V$ en coordonn\'ees de Jacobi est exprim\'ee par:
$$ 
\Delta V=\frac{1}{\sqrt{g}}\left\{\frac{\partial}{\partial \phi}\left(\frac{\sqrt{g}}{g_{11}}.\frac{\partial V}{\partial \phi}\right)+\frac{\partial}{\partial \lambda}\left(\frac{\sqrt{g}}{g_{22}}.\frac{\partial V}{\partial \lambda}\right)+\frac{\partial}{\partial u}\left(\frac{\sqrt{g}}{g_{33}}.\frac{\partial V}{\partial u}\right)\right\} $$
o\`u $g$ est le d\'eterminant de la matrice $G$, donner l'expression de $\Delta V$. 

4. Calculer $\Delta V$ sachant que $V$ est donn\'ee par:
$$ V(\phi,u)=\frac{GM}{\epsilon}Arctg\frac{\epsilon}{u}+ \frac{1}{3}a^2\omega^2\frac{q}{q_0}\left(1-\frac{3}{2}cos^2\phi\right)+\frac{1}{2}\omega^2(u^2+\epsilon^2)cos^2\phi $$
avec $G,M$ et $\omega$ des constantes et:
\ba
q=q(u)=\frac{1}{2}\left[\left(1+3\frac{u^2}{\epsilon^2}\right)Arctg\frac{\epsilon}{u}-3\frac{u}{\epsilon}\right] \nonumber \\
q_0=q(u=b)=\frac{1}{2}\left[\left(1+3\frac{b^2}{\epsilon^2}\right)Arctg\frac{\epsilon}{b}-3\frac{b}{\epsilon}\right] \nonumber
\ea

\newpage
\section{	\textsc{Les R\'eductions des Distances}}

\textbf{Exercice $n^{\circ}$1:}
 On a mesur\'e une distance suivant la pente $D_P = 20130.858\, m $ entre deux points $A$ et $B$ avec  $H_A = 235.07\, m,\, H_B = 507.75\, m$, on prendra comme rayon terrestre  $R = 6378\, km$.

1.	Calculer la distance suivant l'ellipsoide :

-	en utilisant les diff\'erentes corrections,

-	en utilisant la formule rigoureuse.

2.	En prenant la valeur de la formule rigoureuse et sachant que le module lin\'eaire $m$ vaut $0.999\,850\,371$, calculer la distance r\'eduite au plan de la repr\'esentation plane utilis\'ee.
\\

\textbf{Exercice $n^{\circ}$2:} 
Entre 2 points $A$ ( $H_A = 128.26\, m$ ) et  $B$  ( $H_B = 231.84\, m$), la distance $D_P$ suivant la pente est \'egale \`a $15\,498.823\, m$. Soit  $D_0$ la distance corde au niveau de la surface de r\'ef\'erence. L'angle de site observ\'e en $A$ en direction de $B$ est  $i = 0.3523\, gr$.

1.	Calculer la valeur de $D_0$ en utilisant la formule rigoureuse.

2.	Calculer $D_0$ par les corrections.

3.	En adoptant la moyenne des deux m\'ethodes, calculer la distance  $D_e$ r\'eduite \`a la surface de r\'ef\'erence.

4.	Le module lin\'eaire de la repr\'esentation plane Lambert Sud utilis\'ee est de $0.999\,648\,744$. calculer alors la distance $D_r$ r\'eduite au plan de la repr\'esentation.
\\

\textbf{Exercice $n^{\circ}$3:} 
On a mesur\'e une distance suivant la pente entre les points $A\, ( H_A =  1\,319.79\, m)$ et  $B\,( H_B = 1\,025.34\, m)$ avec $D_P = 16\,483.873\, m$.

1. Calculer la distance $D_e$ distance r\'eduite \`a l'ellipsoide de r\'ef\'erence par la formule rigoureuse, on prendra le rayon de la Terre $R= 6378\, km$.

2. Calculer la distance  $D_r$  r\'eduite \`a la repr\'esentation plane Lambert si l'alt\'eration lin\'eaire de la zone est de $- 14\,cm/km$.

\newpage
\section{	\textsc{Les Repr\'esentations Planes}}
\textbf{Probl\`eme $n^{\circ}$1:}
 Soit $\BbS^2$ la sph\`ere de rayon $R$, au point $P(\varphi, \lambda )$ on lui fait correspondre le point $p(X,Y)$ du plan $OXY$ par la repr\'esentation plane suivante d\'efinie par les formules :
	\[                         
	   p (X,Y)=\left\{\begin{array}{ll}
	        X = 2R.tg(\ds \frac{ \pi}{4} -\frac{\varphi}{2}).sin\lambda \\ 
	                        Y = - 2R.tg(\ds \frac{ \pi}{4} -\frac{\varphi}{2}).cos\lambda 
	                        \end{array}\right.
\]
1.	Montrer que l'image d'un m\'eridien ($\lambda$  = constante ) est une droite dont on donne l'\'equation.

2.	Montrer que l'image d'un parall\`ele ($\varphi$  = constante ) est un cercle dont on pr\'ecise l'\'equation.

3.	En utilisant le lemme de Tissot, d\'eterminer les directions principales.

4.	Soit $dS$ la longueur infinit\'esimale correspondante sur le plan, calculer $dS$.

5.  Sachant que sur la sph\`ere $ds^2 = R^2d\varphi^2 +R^2cos^2\varphi .d\lambda ^2$, calculer le module lin\'eaire $m$.

6.	En d\'eduire le module lin\'eaire $m_1$ le long du m\'eridien.

7.	En d\'eduire le module lin\'eaire $m_2$ le long d'un parall\`ele.

8.	Comparer $m_1$  et $m_2$. Conclure sur la conformit\'e ou la non conformit\'e de la repr\'esentation plane.
\\

\textbf{Probl\`eme $n^{\circ}$2:}
 Soit $\Sigma$ la sph\`ere de rayon $R$, au point $P(\varphi, \lambda )$ on lui fait correspondre le point $p(X,Y)$ du plan $OXY$ par la repr\'esentation plane suivante d\'efinie par les formules :
	\[                         
	   p (X,Y)=\left\{\begin{array}{ll}
	        X = R.\lambda  \\ 
	                        Y = R.Logtg(\ds \frac{ \pi}{4}+\frac{\varphi}{2}) 
	                        \end{array}\right.
\]
o\`u   $Log$ d\'esigne le logarithme n\'ep\'erien.

1.	Quelles sont les images des m\'eridiens  ($\lambda$ = constante)  et des  parall\`eles ($\varphi$ = constante). 

2.	Soit $dS$ la longueur infinit\'esimale correspondante sur le plan, calculer $dS$ en fonction de $\varphi$ et  de $\lambda$ et  calculer  le module lin\'eaire $m$.

3.	En d\'eduire les modules lin\'eaires $m_1$ le long du m\'eridien et  $m_2$ le long du parall\`ele.

4.	Comparer $m_1$ et $m_2$  et conclure sur la conformit\'e ou la non conformit\'e de la repr\'esentation plane.

5.	On suppose que $P$ d\'ecrit sur la surface $\Sigma$ une courbe  $(\gamma)$ telle que $\varphi$ et $\lambda$ sont li\'ees par la relation :  $tg\varphi = sin\lambda$. Pour $\varphi = 0\, gr,\, 2\, gr,\, 4\, gr,\, 6\, gr,\, 8\, gr$ et $10\, gr$,  dresser un tableau  donnant  les  valeurs  de  $\lambda$ correspondantes.

6.	Sachant que  $R= 1000\, m$, calculer les coordonn\'ees $(X,Y)$ de la repr\'esentation  plane donn\'ee ci-dessus  pour  les  valeurs de $\varphi$  et $\lambda$  de  la  question 5.

7.	Rapporter \`a l'\'echelle 1/100 sur le plan $OXY$, les positions $(X, Y)$ des points. Que pensez-vous de l'image de la courbe $(\gamma)$. 
\\

\textbf{Probl\`eme $n^{\circ}$3:} 
Sur une sph\`ere de rayon unit\'e, mod\`ele de la terre, on d\'esigne :

-	par  $p$  le p\^ole nord,

-	par $(C)$  un grand cercle qui coupe l'\'equateur au point  $i$  de longitude nulle,

-	par $q$ le p\^ole de ce grand cercle, de latitude $\varphi_0$ positive,

-	par $\omega$  et $h$ respectivement les points d'intersection de $(C)$  et du m\'eridien de $q$  et du grand cercle issu de $q$, passant par le point  $a (\varphi,\lambda)$. 

On pose :    $   \omega h = x,\quad       ha = y$

1.	$q$ est le pivot d'une repr\'esentation cylindrique conforme oblique tangente, dont $(C)$  est le ''pseudo-\'equateur''. Le plan est rapport\'e aux axes   $\Omega X,\Omega  Y$ images respectives de  $(C)$   et du grand cercle   $\omega pq$. Exprimer en fonction de $\varphi, \lambda$    et $\varphi_0$ les coordonn\'ees $X,Y$ du point $A$ image de a (v\'erifier que pour  $\varphi_0 = 0$, on retrouve les expressions  de $X,Y$ d'une repr\'esentation transverse).

2.	Montrer que l'\'equation de l'image plane du parall\`ele de latitude $\varphi_0$  peut s'\'ecrire :                $$e^Y cosX=tg\varphi_0$$                   

Indications : $b$ d\'esignant un point de latitude $\varphi_0$, le triangle $pqb$ est isoc\`ele,  d\'ecomposer ce triangle en deux triangles rectangles \'egaux. Etudier qualitativement les images des autres parall\`eles. 

3.	Montrer que l'image plan de l'\'equateur a pour \'equation :   
                                               $$ cosX + tg\varphi_0.shY = 0$$

Ecrire d'une mani\`ere analogue, l'\'equation de l'image du m\'eridien $\lambda = 0$.

4.	Exprimer le gisement du m\'eridien en fonction de $\varphi,\lambda$   et $\varphi_0$. D\'eterminer la valeur du module lin\'eaire, en particulier en $p$, en un point de l'\'equateur, en un point du m\'eridien origine. 
\\

\textbf{Probl\`eme $n^{\circ}$4:} 
Etude de la repr\'esentation conforme d'une sph\`ere de rayon unit\'e dite repr\'esentation de Littrow\footnote{En hommage \`a  Joseph Johann Littrow (1781-1840) astronome autrichien.} d\'efinie par :
  $$ Z = sinz $$    
	avec  $z =\lambda   + iL$  et  $Z = X + iY$.

1. Pr\'eciser le canevas, les images des m\'eridiens et celle de l'\'equateur.

2. V\'erifier que les points $f$ et $f'$ $(\varphi = 0, \lambda   =\pm \pi/2)$ sont des points singuliers.

3. Etudier les images plans des cercles de diam\`etre $ff'$ et des petits cercles orthogonaux.

4. Soit $s$ le point $(\varphi  =\varphi_0, \lambda   =  0 )$. On appelle segment capable sph\'erique l'ensemble des points $b$ tels que l'angle $ \widehat{bp,bs} = \alpha $. Quelle est l'image plane de cette courbe dans cette repr\'esentation plane.
\\

\textbf{Probl\`eme $n^{\circ}$5:}
 Soit l'application $ F(u,v):\BbR^2\longrightarrow \BbR^3 \verb|\|(0,0,1)$ d\'efinie par:
$$ OM(u,v)=F(u,v)\left\{\begin{array}{l}
x=\ds \frac{2u}{u^2+v^2+1} \\
\\
y=\ds \frac{2v}{u^2+v^2+1}\\
\\
z=\ds \frac{u^2+v^2-1}{u^2+v^2+1}
\end{array}\right. $$ 

1. Calculer la forme fondamentale $ds^2$.

2. Montrer que $OM(u,v)$ appartient \`a la sph\`ere $\BbS^2$ d'\'equation $x^2+y^2+z^2=1$.

3. Calculer $u,v$ en fonction de $x,y$ et $z$.

4. Soit le point $N(0,0,1)$ de  $\BbS^2$, calculer les coordonn\'ees $(X,Y)$ du point $p$ intersection de la droite $NM$ avec la plan $z=0$ en fonction de $x,y$ et $z$. 

5. Soit $\sigma$ l'application $\BbR^3 \verb|\|(0,0,1)\longrightarrow \BbR^2:(x,y,z)\longrightarrow (X,Y)=(X(x,y,z),Y(x,y,z))$. Montrer que $(\sigma \circ F)(u,v)=\sigma(F(u,v))=(u,v)$. En d\'eduire que $F=\sigma^{-1}$.

6. Trouver le rapport de ce probl\`eme avec le probl\`eme \ref{prob1}.
\\

\textbf{Probl\`eme $n^{\circ}$6:} 
Soit un ellipsoide de r\'evolution $E(a,e)$ avec $a$ et $e$ respectivement le demi-grand axe de l'ellipsoide de r\'evolution et $e$ la premi\`ere excentricit\'e. Soit $\BbS^2$ une sph\`ere de rayon $R$. On veut \'etudier le passage suivant: 
$$	p(\varphi,\lambda) \,\,\mbox{de l'ellipsoide}\,\, E \Rightarrow \,P(\psi,\Lambda)\,\mbox{de la sph\`ere }\, \,\BbS^2 $$
1. Exprimer $m$ le module lin\'eaire de cette repr\'esentation.

2. On pose:$$ z=\m L+i\lambda,\quad Z=L+i\Lambda $$
$\m L$ est la latitude isom\'etrique de l'ellipsoide de r\'evolution et $L$ la latitude de Mercator. Une transformation conforme entre $E$ et $\BbS^2$ est donn\'ee par $Z=f(z)$ o\`u $f$ est une fonction analytique. On propose le cas le plus simple \`a savoir:
\ba
Z=\alpha z+\beta \nonumber \\
avec \,\,\left\{\begin{array}{l}
\alpha=c_1+ic_2  \\
\beta=b_1+ib_2
\end{array} \right.\nonumber
\ea
les $c_1,c_2,b_1,b_2$ sont des constantes r\'eelles. Donner les expressions de $L$ et $\Lambda$ en fonction de $\m L$ et $\lambda$.

3. On veut que rep\'esentation transforme les m\'eridiens et les parall\`eles de l'ellipsoide respectivement en m\'eridiens et parall\`eles de la sph\`ere et que l'image du m\'eridien origine $\lambda=0$ soit le m\'eridien origine de la sph\`ere $\Lambda=0$. Montrer que $c_2=b_2=0$ et $ L=c_1\m L+b_1, \quad 
\Lambda= c_1\lambda$.

4. Pour avoir la m\^eme orientation en longitude, on prendra $c_1>0$. On cherchera la transformation \`a d\'eformation minimale autour d'un parall\`ele $\varphi=\varphi_0$ tel que le parall\`ele $\varphi=\varphi_0$ est autom\'ecoique et le module lin\'eaire $m$ est stationnaire pour $\varphi=\varphi_0$, c'est-\`a-dire $m(\varphi_0)=1$ et $\ds \left(\frac{dm}{d\varphi }\right)\biggr |_{\varphi=\varphi_0}=0$, en plus on consid\`ere aussi la condition: $$\ds \left(\frac{d^2m}{d\varphi ^2}\right)\biggr |_{\varphi=\varphi_0}=0 $$
Pour faciliter les notations, on prendra $b=b_1, c=c_1$. Montrer que la relation liant $\varphi_0$ et $\Psi_0$ est:
$$  tg\psi_0=tg\varphi_0\sqrt{\frac{1-e^2}{1-e^2sin^2\varphi_0}} $$
5. D\'eterminer les constantes $b,c$ et $R$ en fonction de $\varphi_0$ et $\Psi_0$ telles que les conditions ci-dessus soient v\'erifi\'ees.

6. Montrer que l'expression du d\'eveloppement limit\'e de $m(\varphi)$ de part et d'autre du parall\`ele $\varphi_0$ est donn\'ee par:
$$m(\varphi)=1-\frac{2e^2(1-e^2)sin\varphi_0cos\varphi_0}{3(1-e^2sin^2\varphi_0)^2}(\varphi -\varphi_0)^3+o((\varphi -\varphi_0)^4)$$
7. On fait intervenir la deuxi\`eme excentricite $e'$, Montrer que $m(\varphi)$ s'\'ecrit:
$$m(\varphi)=1-\frac{2e'^2sin\varphi_0cos\varphi_0}{3(1+e'^2cos^2\varphi_0)^2}(\varphi -\varphi_0)^3+o((\varphi -\varphi_0)^4)$$
\\

\textbf{Probl\`eme $n^{\circ}$7:} 
Soit $\m E(a,b)$ un ellipsoide de r\'ef\'erence de param\`etres $a$ et $e$ respectivement le demi-grand axe et la premi\`ere excentricit\'e. On consid\`ere une repr\'esentation plane $\m P$ de $\m E$ vers le plan $(O,X,Y)$. On pose:
$$\begin{array}{l}
z=\lambda+i\m L\\
Z=X+iY=Z(z)
\end{array}$$
 avec $\m L$ la latitude isoparam\'etrique.

1. Ecrire les expressions du carr\'e des \'el\'ements infinit\'esimaux de longueur sur l'ellipsoide et le plan. En d\'eduire le module lin\'eaire $m$.

2. On pose $\ds \zeta=\frac{\partial Z}{\partial z}$. Si $\gamma$ est le gisement de l'image du m\'eridien passant par le point $z=(\lambda,\m L)$, montrer que $\ds  arg(\zeta)=\frac{\pi}{2}-\gamma$.

3. On cherche une repr\'esentation plane du type $Z=\alpha +\beta z+\varpi z^2$ o\`u $\alpha,\beta$ et $\varpi$ des constantes complexes. On impose les conditions suivantes:

- pour $z=0$, $Z=0$,

- l'axe des $Y$ coincide avec le m\'eridien \`a l'origine.

Montrer que $\m R e (\beta)=0$.

4. En d\'eduire que $Z$ s'\'ecrit:
$$ Z=i\beta_1z+(\varpi_1+i\varpi_2)z^2$$
avec $\beta_1,\varpi_1,\varpi_2$ sont des r\'eels.

\newpage
\section{	\textsc{La Repr\'esentation Lambert}}
\textbf{Exercice $n^{\circ}$1:}
 En un point $A$ de coordonn\'ees g\'eod\'esiques $\varphi = 40.9193\, gr$ et $\lambda = 11.9656\, gr$ \`a l'Est de Greenwich, on vise un point $B$.
 
1.	Dans quelle zone de Lambert Tunisie se trouve le point $A$ ? Calculer ses coordonn\'ees planes $(X,Y)$.
 
2.	L'azimut g\'eod\'esique de la direction $AB$ est $Azg = 55.7631\, gr$. Sachant que $Dv = 1.52\,dmgr$, calculer $G$  le gisement de la direction $AB$.

3.	La distance $AB$ r\'eduite \`a l'ellipsoide de r\'ef\'erence est $D_e = 5421.32\, m$. Sachant que l'alt\'eration lin\'eaire dans la r\'egion des points $A$ et $B$ vaut $- 9\, cm/km$, calculer la distance $AB$ r\'eduite au plan.
\\

\textbf{Exercice $n^{\circ}$2:}
 D'apr\`es les coordonn\'ees de deux points $A$ et $B$ vous trouvez la distance $AB= 5427.380\, m$. Sachant que :

a -	l'alt\'eration lin\'eaire de la repr\'esentation dans la r\'egion de $AB$ vaut $+8.10^{-5}$,

b -	les altitudes des points $A$ et $B$ sont : $H_A = 1000.00\,m$ et $H_B = 1200.00\, m$. Calculer la distance suivant la pente $D_P$ entre les points $A$ et $B$ mat\'erialis\'es sur le terrain.
\\

\textbf{Probl\`eme $n^{\circ}$1:}
 On a mesur\'e une distance suivant la pente $D_P = 20\,130.858\, m$ entre deux points $A$ et $B$ avec  $H_A = 235.07\, m,\, H_B = 507.75\, m$, on prendra comme rayon terrestre  $R = 6378\, km$.

1.	Calculer la distance $D_e$ suivant l'ellipsoide en utilisant la formule rigoureuse.

2.	Sachant que le module lin\'eaire $m$ vaut $0.999\,850\,371$, calculer la distance $D_r$ r\'eduite au plan de la repr\'esentation plane utilis\'ee.

3.	Les coordonn\'ees g\'eod\'esiques du point $A$ sont : $\varphi = 10.7245\,3\, gr,\,\lambda = 41.4490\,3\, gr$. Par des observations astronomiques, on a d\'etermin\'e les coordonn\'ees astronomiques $\varphi_a = 10.7257\,4\, gr$   et $\lambda_a = 41.4505\,2\, gr$  du point $A$ et l'azimut astronomique de la direction $AB$  soit   $Aza = 89.6849\,9\, gr$.	Transformer l'azimut astronomique de la direction $AB$ en azimut g\'eod\'esique en utilisant l'\'equation de Laplace donn\'ee par :$$ Azg = Aza + (\lambda -\lambda_a ).sin\varphi $$
4.	Calculer le gisement $G$ de la direction $AB$ sachant que $Dv = 0.0018\,8\, gr$.

5.	Les coordonn\'ees Lambert Nord Tunisie de $A$ sont $X = 478\,022.43\, m$ et $Y = 444\,702.22\, m$. D\'eterminer alors les coordonn\'ees de $B$.

6.	Calculer l'azimut de $B$ vers $A$ sachant qu'on n\'eglige la correction de la corde de la direction $BA$ et que $\lambda_B =  10.9288\,4 \,gr$. 
\\

\textbf{Probl\`eme $n^{\circ}$2:}
 On a mesur\'e une distance suivant la pente entre les points $A ( H_A =  1319.79 \,m)$ et  $B      ( H_B = 1025.34\, m)$ avec $D_P = 16\,483.873\, m$.

1.	Calculer la distance $D_e$ distance r\'eduite \`a l'ellipsoide de r\'ef\'erence par la formule rigoureuse, on prendra le rayon de la Terre $R= 6378\, km$.

2.	Calculer la distance  $D_r$  r\'eduite \`a la repr\'esentation plane Lambert si l'alt\'eration lin\'eaire de la zone est de $- 14\, cm/km$.

3.	La direction $AB$  a un azimut g\'eod\'esique $Azg = 297.5622\,5\, gr$. Donner l'expression du gisement $G$ de $AB$ en fonction de $Azg,\gamma$ la convergence des m\'eridiens et $Dv$ la correction de la corde, sachant que la repr\'esentation plane utilis\'ee est le Lambert Sud Tunisie et que le point $A$ est au nord du  parall\`ele origine.

4.	On donne $Dv = - 13.7\, dmgr$ et $\lambda    = 9.3474\, 734\, gr$ la longitude de $A$, calculer $G$.

5.	En d\'eduire les coordonn\'ees  $(X_B, Y_B )$ de $B$  si  $X_A = 363\,044.79\, m$  et   $Y_A  = 407\,020.09\, m$.

6.	D\'eterminer les coordonn\'ees g\'eographiques $(\varphi,\lambda)$ de $B$. 

On rappelle que: $a = 6\,378\,249.20\, m$ et $e^2 = 0.006\,803\,487\,7$.

\newpage
\section{	\textsc{La Repr\'esentation UTM}}
\textbf{Exercice $n^{\circ}$1:}
Dans cet exercice, on voudrait justifier l'arr\^et \`a l'ordre 8 de l'expression de $Y(UTM)$ en fonction de $\lambda$. On donne: $\varphi=40.00 \,gr$ et $ a = 6\,378\,249.20\, m,\,\,\,e^2 = 0.006\,803\,4877$.

1. Calculer num\'eriquement $e'^2,\eta^2,t^2=tg^2\varphi$ et $N(\varphi)$.

2. Calculer num\'eriquement le coefficient $a_8$ de (\ref{tabc}).

3. On donne $\lambda= 1.235\,46\,gr$, calculer $a_8\lambda^8$ et conclure.
\\
 
\textbf{Probl\`eme $n^{\circ}$1:}
 Soit le point $A$ de coordonn\'ees g\'eod\'esiques: $\varphi = 40.9193\, gr$ et $\lambda= 11.9656\, gr$ \`a l'Est de Greenwich. On consid\`ere la repr\'esentation plane UTM tronqu\'ee suivante, de m\'eridien central $\lambda_0 = 9^{\circ}$ d\'efinie par les formules :
	\[  \left\{\begin{array}{ll}
	      X =  a_1.(\lambda- \lambda_0)+ a_3. (\lambda- \lambda_0)^3 \\
	        Y =   g(\varphi) + a_2.(\lambda- \lambda_0)^2  
	        \end{array} \right.
\]
o\`u $\varphi,\,\lambda$ et $\lambda_0$ sont exprim\'ees en $rd$, avec: 
	\[  a_1= N(\varphi).cos\varphi \]
	\[a_2 =  \frac{a_1}{2}. sin\varphi \]
	\[a_3 =  \frac{a_1cos^2\varphi}{6}(1- tg^2\varphi + e'^2.cos^2\varphi)\]
	\[N(\varphi) =  \frac{a}{\sqrt{1 - e^2.sin^2\varphi}}\]
	\[g(\varphi) = a(1 - e^2)( 1.0051353.\varphi -  0.0025731sin2\varphi )\]
	 \[ a = 6\,378\,249.20\, m,\,\,\,e^2 = 0.006\,803\,4877,\,\,\, \displaystyle e'^2=\frac{e^2}{1-e^2}   
\]

1.	Montrer que les coordonn\'ees du point $A$ sont: $X = 157\,833.48\, m\,,Y = 4\,078\,512.97\, m$, on justifie les r\'esultats.

2.	Soit le point $B$ de coordonn\'ees $(X = 160\,595.98\,m ; Y = 4\,078\,564.53\, m)$. Sachant que $B$ est situ\'e sur le m\^eme parall\`ele que $A$, calculer la longitude $\lambda'$ de $B$.

3.	Calculer le gisement $G$ et la distance $AB$.

4.	Sachant que la convergence des m\'eridiens $\gamma$  est donn\'ee par $ tg\gamma = (\lambda - \lambda_0)sin\varphi$   et qu'on n\'eglige le $Dv$, calculer l'azimut de la direction $AB$.

5.	Calculer l'azimut de $B$ vers $A$ en n\'egligeant le $Dv$ de $B$ vers $A$.

6.	En calculant les coordonn\'ees UTM de $A$ et $B$, on trouve respectivement  $X_A=657\,770.34\, m,\, Y_A = 4\,076\,891.20\,m;\,  X_B = 660\,531.74\, m,\, Y_B = 4\,076\,942.76\, m$. Calculer la distance $AB$ par les coordonn\'ees UTM. En d\'eduire l'erreur relative sur la distance en utilisant les coordonn\'ees de l'UTM tronqu\'ee.

      \newpage                
\section{	\textsc{Les Transformations de passage entre les Syst\`emes G\'eod\'esiques}}

\textbf{Probl\`eme $n^{\circ}$1:}
 Soient les trois tableaux ci-dessous des coordonn\'ees 3D respectivement dans les syst\`emes $S1$ et $S2$ et \`a transformer dans le syst\`eme $S2$:
\begin{center}
		\[\begin{array}{cccc}
\hline
\hline
 Nom &$X(m)$& $Y(m)$& $Z(m)$\\
\hline
\hline
1& 4\,300\,244.860& 1\,062\,094.681& 4\,574\,775.629\\
2& 4\,277\,737.502& 1\,115\,558.251 &4\,582\,961.996\\
3& 4\,276\,816.431 &1\,081\,197.897 &4\,591\,886.356\\
4 &4\,315\,183.431& 1\,135\,854.241& 4\,542\,857.520\\
5& 4\,285\,934.717& 1\,110\,917.314& 4\,576\,361.689\\
6 &4\,217\,271.349& 1\,193\,915.699& 4\,618\,635.464\\
7& 4\,292\,630.700& 1\,079\,310.256& 4\,579\,117.105\\
 \hline
\hline
\end{array}\]
	\vspace{1cm}
		\[\begin{array}{cccc}
\hline
\hline
 Nom &$X(m)$& $Y(m)$& $Z(m)$\\
\hline
\hline
1 &4\,300\,245.018& 1\,062\,094.592 &4\,574\,775.510\\
2 &4\,277\,737.661& 1\,115\,558.164 &4\,582\,961.878\\
3 &4\,276\,816.590& 1\,081\,197.809& 4\,591\,886.238\\
4 &4\,315\,183.590 &1\,135\,854.153& 4\,542\,857.402\\
5 &4\,285\,934.876 &1\,110\,917.227& 4\,576\,361.571\\
6 &4\,217\,271.512 &1\,193\,915.612& 4\,618\,635.348\\
7 &4\,292\,630.858 &1\,079\,310.168& 4\,579\,116.986\\
\hline
\hline
\end{array}\]
		\[\begin{array}{cccc}
\hline
\hline
 Nom &$X(m)$& $Y(m)$& $Z(m)$\\
\hline
\hline
A &4\,351\,694.594& 1\,056\,274.819&4\,526\,994.706\\
B &4\,319\,956.455& 1\,095\,408.043& 4\,548\,544.867\\
C &4\,303\,467.472 &1\,110\,727.257 &4\,560\,823.460\\
D &4\,202\,413.995 &1\,221\,146.648 &4\,625\,014.614\\
\hline
\hline
\end{array}\]
\end{center}
1. D\'eterminer les param\`etres du mod\`ele de Bur$\breve{s}$a-Wolf \`a 7 param\`etres.

2. Calculer les coordonn\'ees 3D des points du troisi\`eme tableau dans le syst\`eme $S2$. 
\
\section{	\textsc{Notions sur le Mouvement d'un Satellite Artificiel de la Terre}}
\textbf{Exercice $n^{\circ}$1:} 
1. Montrer que: $r=a(1-ecosE) $.

2. D\'emontrer \`a partir des formules du cours la relation:
$$ tg\frac{\upsilon}{2}=\sqrt{\frac{1+e}{1-2}}tg\frac{E}{2}$$
Aide: exprimer $tg(\upsilon/2)$ en fonction de $tg\upsilon $.
 \\

\textbf{Exercice $n^{\circ}$2:}
 A partir de l'expression de $X_C$, monter que $X_C$ v\'erifie l'\'equation du mouvement non perturb\'e pour la composante $X$, soit: 
$$ \ddot{X}_C+\frac{\mu}{r^3}X_C=0 $$
\\

\textbf{Probl\`eme $n^{\circ}$1:} 
La Terre est suppos\'ee sph\'erique, homog\`ene de rayon $R = 6\,371\,000\,m$. Le produit de la constante universelle de gravitation
 terrestre $G$ par la masse $M$ de la Terre soit $GM = 3.986\, 005\,10^{14}\,m^3s^{-2}$. Un satellite g\'eod\'esique a une trajectoire telle que son altitude maximale est $1100\,km$ et son altitude minimale $800\,km$.

1. Donner la p\'eriode de ce satellite.

2. Quelle est l'excentricit\'e de sa trajectoire?

3. On mesure la distance du satellite \`a une station au sol de latitude $43^{\circ},5$ et d'altitude nulle, lors du passage du satellite \`a la verticale de la station, soit $D = 812\,000\,m$.

a - Quelle est l'anomalie vraie du satellite \`a cet instant, sachant qu'il vient de passer au p\'erig\'ee.

b - Combien de temps s'est \'ecoul\'e depuis le passage au p\'erig\'ee?
\\

\textbf{Probl\`eme $n^{\circ}$2:} 
Une com\`ete d\'ecrit autour du Soleil une ellipse d'excentricit\'e $e$ de demi-grand axe $a$ et de demi-petit axe $b$ o\`u le Soleil occupe un des foyers. L'\'equation de l'orbite de la com\`ete en coordonn\'ees polaires est donn\'ee par:
$$ r=\frac{a(1-e^2)}{1+ecos\upsilon}$$ avec $r$ la distance Soleil- com\`ete.

1. D\'eterminer les distances $r_A$ et $r_P$ lorsque la com\`ete est \`a l'apog\'ee et au p\'erig\'ee en fonction de $a$ et $e$.

2. La com\`ete de Halley a une orbite fortement excentrique : son apog\'ee est \`a 0.53 $UA$  du Soleil et sa p\'erig\'ee est \`a 35.1 $UA$. Calculer $e$.

3. En utilisant la loi des aires et la troisi\`eme loi de Kepler, montrer que la constante des aires $C$ est exprim\'ee par:
$$C^2=\frac{b^2}{a}G.M$$
o\`u $G, M$ d\'esignent respectivement la constante de la gravitation universelle et la masse du Soleil.

4. On pose : $\ds u=\frac{1}{r}$. Donner l'expression du carr\'e de la vitesse $v^2$ de la com\`ete en fonction de $u$ et $\ds \frac{du}{d\upsilon}$. Montrer que $v^2$ peut s'\'ecrire sous la forme:
$$ v^2=G.M\left(\frac{2}{r}-\frac{1}{a}\right)$$

5. D\'eterminer l'expression du rapport des vitesses \`a l'apog\'ee et au p\'erig\'ee $\ds \frac{v_A}{v_P}$ en fonction de $e$.

6. Calculer num\'eriquement ce rapport pour le cas de la com\`ete de Halley.

On donne:

- 1 $UA=149\, 597\,870 \,km$, 

- $G=6.672\times 10^{-11}\,m^3.kg^{-1}.s^{-2}$,

- $M=1.9891\times 10^{30}\,kg$.

\part{\textsc{Th\'eorie des Erreurs}}
\chapter{Exercices et Probl\`emes}
\section{\textsc{Exercices et Probl\`emes}}
\textbf{Exercice $n^{\circ}$1:}
 Soient $U$ un ouvert convexe d'un espace de Banach\footnote{\textbf{Stefan Banach} (1892-1945): math\'ematicien polonais.} r\'eel $E$ c'est-\`a-dire un espace vectoriel norm\'e complet sur $\BbR$ et $f$ une fonction \`a valeurs r\'eelles, diff\'erentiable et convexe dans $U$. Montrer que si $f'(x_0)=0$ en un point $x_0\in U$, alors $f$ a un minimum absolu en $x_0$.\index{\textbf{Banach S.}}
\\

\textbf{Exercice $n^{\circ}$2:}
 Montrer que dans un espace de Banach r\'eel $E$, la fonction $f=\|\,. \|^2$ est strictement convexe, c'est-\`a-dire, $\forall \, \alpha \in ]0,1[,\,\,f(\alpha x+(1-\alpha) y)< \alpha f(x)+(1-\alpha)f(y), $ pour tout couple $(x,y) \in E^2$.

Aide: utiliser l'identit\'e remarquable:
$$\|\alpha x + (1-\alpha )y\|^2 = \alpha \|x\|^2 + (1-\alpha )\|y\|^2-\alpha (1-\alpha )\|x- y\|^2$$
\\

\textbf{Exercice $n^{\circ}$3:}
 On note $F$ une surface de $\BbR^3$ d\'efinie par la repr\'esentation param\'etrique:
$$\textbf{\textit{OM}}=(a_1(u,v),a_2(u,v),a_3(u,v))^T $$
o\`u $u$ et $v$ sont deux param\`etres r\'eels. On se donne un point $P(x,y,z)\in \BbR^3$.

1. Donner une condition g\'eom\'etrique portant sur le plan tangent \`a $F$ au point $M_0(u_0,v_0)$ pour que la diff\'erentielle de la fonction $(u,v)\longrightarrow \varphi(u,v)=\left\|\textbf{\textit{OP}}-\textbf{\textit{OM}}(u,v)\right\|^2$ soit nulle en $M_0(u_0,v_0)$.
\\

\textbf{Exercice $n^{\circ}$4:} 
Soient $U$ un ouvert convexe d'un espace de Banach r\'eel $E$ et $f$ une application diff\'erentielle de $U$ dans $\BbR$.

1. Montrer que $f$ est convexe dans $U$ si et seulement si:
$$ f(x)\geq f(x_0)+f'(x_0)(x-x_0)$$
pour tout couple de points $x,x_0\in U$.

2. On suppose $E=\BbR^n$ et $f$ de classe $C^2$ soit deux fois diff\'erentiable et $f"$ continue; pour $x\in U$, soit $\varphi_x$ la forme quadratique  d\'efinie par :
$$\varphi_x(h)= \sum_{i,j=1}^n\frac{\partial ^2 f}{\partial x_i \partial x_j}(x)h_ih_j,\quad \,\,\,\,h=(h_1,h_2,...,h_n)\in \BbR^n $$
Montrer que $f$ est convexe dans $U$ si et seulement si $\varphi_x$ est positive pour tout $x\in U$ soit $\varphi_x(h)\geq 0$ pour $x\in U$ et $h\in \BbR^n$. 
\\

\textbf{Exercice $n^{\circ}$5:} 
Soit un triangle $ABC$, on observe les angles $\hat{A},\,\hat{B},\,\hat{C}$ et les c\^ot\'es $BC=a,\, AC=b$ et $AB=c$:
	\[
\left\{\begin{array}{l}
\hat{A}=43.7716\,0\,gr\,\quad \sigma_{\hat{A}}=3.1\,dmgr\\
\hat{B}=98.3904\,3\,gr\,\quad \sigma_{\hat{B}}=3.1\,dmgr\\
\hat{C}=57.8385\,8\,gr\,\quad \sigma_{\hat{C}}=3.1\,dmgr\\
a=333.841\,m,\quad \sigma_a=0.005\,m \\
b=525.847\,m,\quad \sigma_b=0.010\,m\\
c=414.815\,m,\quad \sigma_c=0.005\,m 
\end{array}\right. \nonumber
\]
1. Calculer les angles et les c\^ot\'es compens\'es.

2. Calculer les poids de l'angle $\hat{A}$ et du c\^ot\'e $a$.

3. D\'eterminer une estimation du facteur de variance unitaire.
\\

\textbf{Probl\`eme $n^{\circ}$1:} 
Les directions suivantes sont observ\'ees respectivement aux stations $A,B,C$ et $D$ d'un quadrilat\`ere $ABDC$ comme suit:
 \ba
Station\,\,A=\left\{\begin{array}{l}
vers\,\,B: \,\,0.0000\,0\,gr \\
vers\,\,C: 74.1666\,7\,gr  
\end{array}\right. \nonumber \\
Station\,\,B=\left\{\begin{array}{l}
vers\,\,D: \,\,0.0000\,0\,gr \\
vers\,\,C: 82.4608\,0 \,gr \\
vers\,\,A: 170.6253\,1\,gr 
\end{array}\right. \nonumber \\
Station\,\,C=\left\{\begin{array}{l}
vers\,\,A: \,\,0.0000\,0\,gr \\
vers\,\,B: 37.6709\,9\,gr \\
vers\,\,D:85.0830\,2\,gr 
\end{array}\right. \nonumber \\
Station\,\,D=\left\{\begin{array}{l}
vers\,\,C: \,\,0.0000\,0\,gr \\
vers\,\,B: 70.1280\,9\,gr 
\end{array}\right. \nonumber
 \ea
Les observations sont non corr\'el\'ees. l'\'ecart quadratique moyen de ces observations est identique et vaut $\sigma_d=6.2\,dmgr$. 

1. Compenser les directions et calculer leurs poids et celui de l'angle $CBA$.

2. Calculer l'estimateur $s^2$ du facteur de variance unitaire et celui de $\ds \frac{s^2}{\sigma^2}$.

3. Des observations de nivellement ont \'et\'e effectu\'ees sur les lignes $ABC$ et $BCD$. Les diff\'erences d'altitudes observ\'ees sont les suivantes:
\ba
H_A-H_B=0.509\,m \nonumber \\ 
H_B-H_D=1.058\,m \nonumber \\ 
H_A-H_C=3.362\,m \nonumber \\ 
H_D-H_C=1.783\,m \nonumber \\ 
H_B-H_C=2.829\,m \nonumber  
\ea 
Les observations sont non corr\'el\'ees et de pr\'ecision identique. Compenser les observations ci-dessus et calculer un estimateur du facteur de variance unitaire.
\\

\textbf{Probl\`eme $n^{\circ}$2:}
 1. Montrer que dans un cheminement altim\'etrique de pr\'ecision, le poids de l'observation entre deux rep\`eres est inversement proportionnel de leur distance en supposant l'\'egalit\'e des port\'ees et que les observations sont non corr\'el\'ees. 

\begin{figure}
	\centering
		\includegraphics{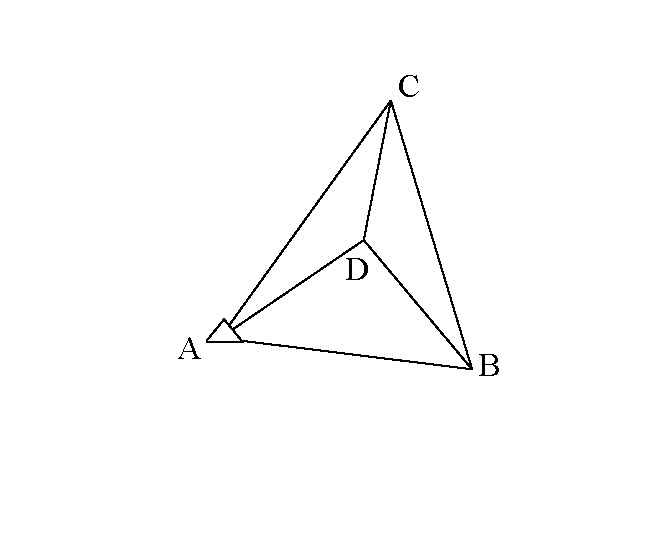}
	\caption{La polygonale observ\'ee}
	\label{fig:triangleniv}
\end{figure}

2. Une polygonale $ABCD$ (voir \textbf{Fig. \ref{fig:triangleniv}}) a \'et\'e observ\'ee par le nivellement de pr\'ecision. L'instrument utilis\'e a une pr\'ecision de $2\,mm$ par $km$. Les observations consid\'er\'ees non corr\'el\'ees sont les suivantes:
\ba
H_C-H_A=1.878\,m, \quad 	AC=6.44\,km \nonumber \\ 
H_D-H_A=3.831\,m, \quad 	AD=3.22\,km \nonumber \\ 
H_D-H_C=1.954\,m, \quad 	CD=3.22\,km \nonumber \\ 
H_B-H_A=0.332\,m, \quad 	AB=6.44\,km \nonumber \\ 
H_D-H_B=3.530\,m, \quad 	BD=3.22\,km \nonumber \\ 
H_C-H_B=1.545\,m, \quad 	BC=6.44\,km \nonumber  
\ea 

L'altitude du rep\`ere $A$ est de $3.048\,m$ et non entach\'ee d'erreurs. Calculer par compensation des observations les altitudes des rep\`eres $B,C$ et $D$ et leurs \'ecarts-types.

3. Calculer l'\'ecart-type de la diff\'erence d'altitudes entre les rep\`eres $C$ et $D$.

4. Donner une estimation de la pr\'ecision par $km$ du nivellement effectu\'e.
\\

\textbf{Probl\`eme $n^{\circ}$3:} 
		On veut \'etalonner un an\'eroide, appareil donnant la pression de l'air, par la formule:
	\[ D=d+\alpha t+\gamma
\]
o\`u $\alpha,\gamma$ sont deux constantes, $t$ la temp\'erature en degr\'es centigrades. Les param\`etres $d$ et $D$ sont lus respectivement de l'an\'eroide et \`a partir d'un barom\`etre en mercure, et exprim\'es en $mm$ $Hg$. 
\\

Pour d\'eterminer $\alpha$ et $\gamma$, des lectures ont \'et\'e prises \`a diff\'erentes temp\'eratures (voir tableau \ref{tab: Tableau des observations}). 
\begin{table}[ht]
	\centering
\begin{center}
\[\begin{array}{lll}
\hline
 \quad  \quad	t		 & d   &      D   \\ \hline
^{\circ}\,\mbox{Centigrade} \,\,      &    mm\,\,Hg \,\,   &              mm\,\,Hg     \\   \hline
 \quad  \quad6.0         &       761.3                 &       762.3              \\   \hline
 \quad  \quad10.0        &      759.1                 &        759.5              \\  \hline 
  \quad  \quad14.0       &    758.4                    &    758.7                  \\ \hline
 \quad  \quad 18.0       &    763.1                    &    763.0                  \\ \hline
\end{array} 
\]
\caption{Table des observations}
	\label{tab: Tableau des observations}
\end{center}
\end{table}
Ces observations sont non corr\'el\'ees. L'\'ecart-type de la lecture de $d$ est de $0.14\,mm\,\,Hg$; $t$ et $D$ sont suppos\'ees sans erreurs.

1. Calculer par la m\'ethode des moindres carr\'es les constantes $\alpha$ et $\gamma$.

2. Estimer le facteur de variance unitaire.

3. D\'eterminer la variance et la covariance de $\alpha$ et $\gamma$.
\\

\textbf{Probl\`eme $n^{\circ}$4:} 
En statistiques, la loi normale est une famille de distributions de probabilit\'es caract\'eris\'ees par la fonction de densit\'e:
$$p(x,\mu, \sigma)=\frac{1}{\sqrt{2\pi}\sigma}e^{-\frac{(x-\mu)^2}{2\sigma^2}}$$
o\`u $\mu$ est la moyenne et $\sigma^2$ la variance. On note par $l(x,\mu,\sigma)=Logp(x,\mu,\sigma)$, soit: $$ l(x,\mu,\sigma)=-Log\sigma-\frac{(x-\mu)^2}{2\sigma^2}$$
Soit $X$ une variable al\'eatoire ayant comme fonction de densit\'e $p(x,\mu,\sigma)$. On rappelle les op\'erateurs suivants esp\'erance math\'ematique ou moyenne et variance:
\ba
E[f(X)]=\int_{-\infty}^{+\infty}f(x)p(x,\mu,\sigma)dx\nonumber \\
V(f(X))=E[(E[f(X)]-f(X))^2]\nonumber
\ea
On donne la formule: $\ds \int_0^{+\infty}e^{-u^2}du=\frac{\sqrt{\pi}}{2}$.

1. Montrer que:
\ba
E(X)=\ds \int^{+\infty}_{-\infty}p(x,\mu,\sigma)dx=\mu \nonumber \\
\sigma^2(X)=Var(X)=Cov(X,X)=\ds \int^{+\infty}_{-\infty}(x-\mu)^2p(x,\mu,\sigma)dx=\sigma^2 \nonumber
\ea 
2. Montrer que:
$$\ds \int_{-\infty}^{+\infty}u^4e^{-u^2}du=\ds \frac{3\sqrt{\pi}}{4}$$

\noindent 3. Calculer $\ds \frac{\partial l}{\partial \mu}\,\,\frac{\partial l}{\partial \sigma}$.
\\

\noindent 4. On pose $\theta=(\mu,\sigma)$. Soit $T_\theta$ l'espace engendr\'e par $(\ds \frac{\partial l}{\partial \mu},\,\,\frac{\partial l}{\partial \sigma})$. On d\'efinit sur $ T_\theta$ l'op\'erateur $<.,.>: T_\theta \times T_\theta\longrightarrow \BbR$ \`a $A,B$ deux variables al\'eatoires $\in T_\theta$:
$$<A,B>=E[A(x)B(x)]$$
Justifier qu'on peut \'ecrire:
$$E[A(x)B(x)]=Cov(A(x),B(x))=E[(E[A(x)]-A(x))(E[B(x)]-B(x))]$$
5. Montrer que $<.,.>$ d\'efinit un produit scalaire sur $T_\theta$.
\\

\noindent 6. On pose: $e_1=\ds \frac{\partial l}{\partial \mu}$ et $e_2=\ds \frac{\partial l}{\partial \sigma}$. On d\'efinit le tenseur m\'etrique sur $T_\theta$ par:
$$g_{ij}=<e_i,e_j>$$
Montrer que la matrice $g=(g_{ij})$ est donn\'ee par: 
$$g=\ds \frac{1}{\sigma^2}\begin{pmatrix}
1 & 0 \\
0 & 2 
\end{pmatrix}$$
et que la premi\`ere forme fondamentale sur $T_\theta$ s'\'ecrit:
$$ ds^2=\ds \frac{1}{\sigma^2}(d\mu^2+2d\sigma^2)$$
\\

\textbf{Probl\`eme $n^{\circ}$5:} 
Soit un triangle de c\^ot\'es $a,b,c$ et d'angles $A,B$ et $C$. On se propose:

- d'estimer $\dot{a},\dot{b}$ et $\dot{c}$,  et les variances de ces d\'eterminations. Les observations sont:
\be
\left\{\begin{array}{l}
   a=96.48\,mm \\
	b=115.50\,mm \\
 	A=63.042 \,gr  \\
  B=99.802\,gr     \\
   C=37.008\,gr    
	\end{array}\right. \lb{AF21}
 \ee
 On choisit ici comme \un{unit\'es normalis\'ees} le d\'ecimillim\`etre $(0.1\,mm)$ pour les mesures de distances, et le d\'ecimilligrade $(0.1\,gr)$ pour les angles. 
\\

On prend les poids \'egaux aux inverses des carr\'es des $emq$ de chaque observation. On donne la matrice des poids $P$:
$$ P=\begin{pmatrix}
0.277 & 0&0&0&0 \\
0&0.160&0&0&0 \\
0&0&1.524&0&0 \\
0&0&0&1.524 &0 \\
0&0&0&0&1.524 
\end{pmatrix}$$
On prendra comme valeurs approch\'ees des inconnues $ a_0=a;\quad b_0=b;\quad c_0=\ds a\frac{sinC}{sinA}$.
\\

1. Ecrire les param\`etres observ\'ees et les valeurs observ\'ees des inconnues dans les nouvelles unit\'es. 

2. Soit $X=(a,b,c)$ le vecteur des inconnues. On adopte le syst\`eme suivant liant les inconnues aux observables:
\be
\left\{\begin{array}{l}
\dot{a}=\dot{a} \\
\dot{b}=\dot{b} \\
Arccos\ds \frac{\dot{b}^2+\dot{c}^2-\dot{a}^2}{2\dot{b}\dot{c}}=\dot{A} \\
\\
 Arccos\ds \frac{\dot{c}^2+\dot{a}^2-\dot{b}^2}{2\dot{c}\dot{a}}=\dot{B} \\
\\
Arccos\ds \frac{\dot{a}^2+\dot{b}^2-\dot{c}^2}{2\dot{a}\dot{b}}=\dot{C} 
\end{array}\right. \lb{BF21-13}
\ee
Ceci \'etant, on posera pour les grandeurs \`a d\'eterminer:
\ba
\dot{a}=a_0+da=a+da \nonumber \\
\dot{b}=b_0+db=b+db \nonumber \\
\dot{c}=c_0+dc \nonumber
\ea
et pour les grandeurs observ\'ees:
\ba
\dot{a}=a+v_a\nonumber \\
\dot{b}=b+v_b\nonumber \\
\dot{A}=A+v_A\nonumber \\
\dot{B}=B+v_B \nonumber\\
\dot{C}=C+v_C \nonumber 
\ea
En linearisant la troisi\`eme \'equation de (\ref{BF21-13}), montrer que l'\'equation d'observation s'\'ecrit:
$$\frac{1}{sinA}\frac{a_0}{b_0c_0}\frac{2000}{\pi}da-\frac{1}{sinA}\frac{a_0^2+b_0^2-c_0^2}{2b_0^2c_0}\frac{2000}{\pi}db-\frac{1}{sinA}\frac{a_0^2+c_0^2-b_0^2}{2b_0c_0^2}\frac{2000}{\pi}dc=-k_A\frac{2000}{\pi}+v_A$$
o\`u :
$$k_A=\frac{b^2_0+c^2_0-a^2_0-2b_0c_0cosA}{2b_0c_0sinA}$$
(\'etant entendu qu'on exprime $v_A$ en $dcgr$).
\\

3. Montrer que le syst\`eme des moindres carr\'es $ AX=L+V$ s'\'ecrit:
$$
\begin{pmatrix}
1. & 0.& 0. \\
0. & 1. & 0. \\
1.00375 & -0.83924 & 0.00143 \\
-1.00571 & 1.20285 & -0.66128 \\
0.00094 & -0.36239 & 0.65918 
\end{pmatrix}.\begin{pmatrix} 
da \\
db \\
dc 
\end{pmatrix}=\begin{pmatrix}
0. \\
0. \\
0.97981 \\
-2.88449 \\
0.42396
\end{pmatrix}+\begin{pmatrix}
v_a \\
v_b \\
v_A \\
v_B \\
v_C 
\end{pmatrix} $$
4. R\'esoudre le syst\`eme pr\'ec\'edent par la m\'ethode des moindres carr\'es et montrer que la matrice normale $N=A^TPA$ est donn\'ee par:
$$ N=\begin{pmatrix} 
3.35605 & -3.13044 & 1.01750 \\
-  & 3.64132 & -1.57937 \\
- & - & 1.32971 
\end{pmatrix}$$
5. Montrer que:
$$ X=\begin{pmatrix}
+0.62971\\
-0.90962\\
 0.94782 
\end{pmatrix}$$
6. D\'eterminer les variances des inconnues $\sigma^2_a,\,\sigma^2_b$ et $\sigma^2_c$.
\\

\textbf{Exercice $n^{\circ}$1:}
 On consid\`ere $(u,v)\in \BbR^2$ et on d\'efinit la fonction par :
$$ f(u,v)=u^4+6uv+1.5v^2+36v+405 $$
1. Chercher les points critiques r\'eels de $f$.

2. Montrer que le point $x^*=(u,v)=(3,-18)$ est un point minimum de $f$.

3. Montrer que le Hessien de $f$ est une matrice d\'efinie positive si $u^2>1$ et ind\'efinie si $u^2< 1$.

4. Montrer que la formule de r\'ecurrence de Newton s'\'ecrit avec $ J=1.5(u_k^2-1)$:
	\[u_{k+1}=\frac{u_{k}^3+9}{J},\quad v_{k+1}=-\frac{2u_{k}^3+18u_k^2}{J}
	\]
	\\

\textbf{Probl\`eme $n^{\circ}$6:} 
Soient le plan $(P)$ et la sph\`ere $(\BbS^2)$ d'\'equations respectivement: $x+y+z=1$ et $x^2+y^2+z^2=1$. On veut chercher le point $M\in (\BbS^2)$  tel que sa distance au plan $(P)$ soit maximale.

1. Montrer que la distance d'un point $M(X,Y,Z)$ au plan $(P)$ est donn\'ee par : 
	\[d=|X+Y+Z-1|/\sqrt{3} 	\]
2. Pour r\'epondre \`a la question pos\'ee ci-dessus, on consid\`ere la fonction: $E(x,y,z,\lambda)=-(x+y+z-1)^2$$-\lambda(x^2+y^2+z^2-1)$. Ecrire le syst\`eme d'\'equations donnant les points critiques de $E$ qu'on note par (1).

3. Montrer que si $\lambda=-1$, on arrive \`a une contradiction. On suppose que $\lambda \neq-1$. Que repr\'esente le cas $\lambda=0$.

4. On suppose que $\lambda \notin \{-3,-1,0\}$. R\'esoudre le syst\`eme (1). Soit le point $M_2$ tel que ses coordonn\'ees sont n\'egatives.

5. Montrer que la matrice hessienne de $E$ pour $M_2$ s'\'ecrit sous la forme:
	\[ H=\begin{pmatrix}
\mu^2 & -2 & -2 \\
-2& \mu^2 & -2 \\
-2 & -2 & \mu^2 
\end{pmatrix}\quad avec\,\,\,\mu=1+\sqrt{3}\]
6. Si on pose $U=(X,Y,Z)^T \in (\BbS^2)$. Montrer que $U^T.H.U=2\left[3+\sqrt{3}-(X+Y+Z)^2\right]$. En d\'eduire que $U^T.H.U>0$ pour tout $U\neq 0 \in (\BbS^2)$. 

7. Montrer que pour le point $M_2$, on obtient un minimum strict de $E$. A-t-on r\'epondu \`a la question du probl\`eme. 
\\

\textbf{Probl\`eme $n^{\circ}$7:}
 Dans le plan affine $\mathcal P$, on a mesur\'e trois distances planes entre un point inconnu $P(X_1,X_2)$ vers trois points connus $P_i(a_i,b_i)_{i=1,3}$ dans trois directions diff\'erentes. On consid\`ere le mod\`ele non lin\'eaire de Gauss-Markov d\'efini par:
$$	\zeta(X)=L-e, \quad e \in \mathcal{N}(0,\Gamma) $$
avec:

- $L$: le vecteur des observations $(3\times1)=(L_1,L_2,L_3)^T$;

- $X$: le vecteur des inconnues $(2\times1)=(X_1,X_2)^T$;

- $e$: le vecteur des erreurs $(3\times1)=(e_1,e_2,e_3)^T$ suit la loi normale $\mathcal{N}(0,\Gamma)$ avec $E(e)=0$ et $\Gamma=E(ee^T)$ la matrice de dispersion ou variance, on prendra $\Gamma=\sigma^2_0.P^{-1}$, $P$ est la matrice des poids \'egale \`a la matrice unit\'e $I_3$,  $\sigma_0$ une constante positive;

- $\zeta$: est une fonction donn\'ee injective d'un ouvert $U\subset \BbR^2 \rightarrow \BbR^3$ d\'efinie par:
$$\zeta(X)=\zeta(X_1,X_2)=\begin{pmatrix} 
\frac{1}{2}\left[(X_1-a_1)^2+(X_2-b_1)^2\right] \\
\frac{1}{2}\left[(X_1-a_2)^2+(X_2-b_2)^2\right] \\
\frac{1}{2}\left[(X_1-a_3)^2+(X_2-b_3)^2\right] 
\end{pmatrix} $$
On prendra comme composante $L_i$ du vecteur  observation la quantit\'e $L_i=\ds \frac{D^2_{i\,\mbox{observ\'ee}}}{2}$

1. Montrer que les vecteurs $\displaystyle \frac{\partial \zeta}{\partial X_1},\frac{\partial \zeta}{\partial X_2}$ sont lin\'eairement ind\'ependants en chaque point $X \in U$.

2. Montrer que les fonctions $\displaystyle \frac{\partial^2 \zeta}{\partial X_i \partial X_j}$ sont continues sur $U$ pour $i,j \in \left\{1,2\right\}$.

3. Posons: $ J=\|L-\zeta(X)\|^2$

Calculer les coefficients de la matrice $(\displaystyle \frac{\partial^2 J}{\partial X_i \partial X_j} ),i,j\in \left\{1,2\right\} $.

4. Soit la matrice carr\'ee d\'efinie par:
$$	g(X)=(g_{ij}) \quad \mbox{avec}\,g_{ij}=<\frac{\partial \zeta(X)}{\partial X_i},\frac{\partial \zeta(X)}{\partial X_j}> \quad \left\{
\begin{array}{ll}
	i=1,2 \\
	j=1,2
\end{array} \right.$$
Calculer les coefficients $g_{ij}$.

5. Introduisons la matrice $B$ d\'efinie par:
$$	B(X,L)=(B_{ij}) \quad \mbox{avec}\,\,\,B_{ij}= g_{ij}-<L-\zeta(X),\frac{\partial^2 \zeta}{\partial X_i\partial X_j}> \quad \left\{
\begin{array}{ll}
	i=1,2 \\
	j=1,2
\end{array} \right. $$
Calculer les \'el\'ements de la matrice $B$ et montrer qu'elle est d\'efinie positive.

\addtocontents{toc}{\protect\vspace*{\baselineskip}}



\clearpage


\end{document}